\documentclass{elsart}
\usepackage{amsfonts}
\usepackage{mathrsfs}
\usepackage{amssymb,amsmath}
\usepackage{graphics}
\usepackage{color}
\usepackage{graphicx}
\numberwithin{thm}{section}\allowdisplaybreaks
\numberwithin{equation}{section}\allowdisplaybreaks

\begin{document}

\begin{frontmatter}

\title{Global well posedness and scattering for the elliptic and non-elliptic derivative nonlinear Schr\"odinger
equations with small data}

\author{Wang Baoxiang\corauthref{*}}\ead{wbx@math.pku.edu.cn}
\corauth[*]{Corresponding author.}

\address{LMAM, School of Mathematical
Sciences, Peking University, Beijing 100871, People's Republic of
China}

\date{March 10, 2008}

\begin{abstract}
\rm  We study the Cauchy problem for the generalized elliptic and
non-elliptic derivative nonlinear Schr\"odinger equations,  the
existence of the scattering operators and the global well posedness
of solutions with small data in Besov spaces
$B^s_{2,1}(\mathbb{R}^n)$  and in modulation spaces
$M^s_{2,1}(\mathbb{R}^n)$ are obtained. In one spatial dimension, we
get the sharp well posedness result with small data in critical
homogeneous Besov spaces $\dot B^s_{2,1}$. As a by-product, the
existence of the scattering operators with small data is also shown.
In order to show these results, the global versions of the estimates
for the maximal functions on the elliptic and non-elliptic
Schr\"odinger groups are established.
\end{abstract}


\end{frontmatter}

\section{Introduction}
We consider the Cauchy problem for the generalized derivative
nonlinear Schr\"odinger equation (gNLS)
\begin{align}
{\rm i} u_t + \Delta_\pm u = F(u, \bar{u}, \nabla u, \nabla
\bar{u}), \quad u(0,x)= u_0(x),
 \label{gNLS}
\end{align}
where $u$ is a complex valued function of $(t,x)\in \mathbb{R}
\times \mathbb{R}^n$,
\begin{align}
\Delta_\pm u = \sum^n_{i=1} \varepsilon_i \partial^2_{x_i}, \quad
\varepsilon_i \in \{1,\,  -1\}, \quad i=1,...,n,
 \label{Delta-pm}
\end{align}
$\nabla =(\partial_{x_1},..., \partial_{x_n})$, $F:
\mathbb{C}^{2n+2} \to \mathbb{C}$ is a polynomial,
\begin{align}
F(z) = P(z_1,..., z_{2n+2})= \sum^{}_{m+1\le |\beta|\le M+1} c_\beta
z^\beta, \quad c_\beta \in \mathbb{C},
 \label{poly}
\end{align}
$m, M\in \mathbb{N}$ will be given below.

There is a large literature which is devoted to the study of
\eqref{gNLS}. Roughly speaking, three kinds of methods have been
developed for the local and global well posedness of \eqref{gNLS}.
The first one is the energy method, which is mainly useful to the
elliptic case $\Delta_\pm=\Delta=\partial^2_{x_1}+...+
\partial^2_{x_n}$, see Klainerman \cite{Klai}, Klainerman and
Ponce \cite{Kl-Po}, where the global classical solutions were
obtained for the small Cauchy data with sufficient regularity and
decay at infinity, $F$ is assumed to satisfy an energy structure
condition ${\rm Re} \;\partial F/\partial (\nabla u) =0$. Chihara
\cite{Chih1,Chih2} removed the condition ${\rm Re}\;
\partial F/\partial (\nabla u) =0$ by using the smooth operators and
the commutative estimates between the first order partial
differential operators and ${\rm i} \partial_t+ \Delta$, suitable
decay conditions on the Cauchy data are still required in
\cite{Chih1,Chih2}. Recently, Ozawa and Zhang \cite{Oz-Zh} removed
the assumptions on the decay at infinity of the initial data. They
obtained that if $n\ge 3$, $s>n/2+2$,  $u_0 \in H^s$ is small
enough, $F$ is a smooth function vanishing of the third order at
origin with ${\rm Re}\;
\partial F/ \partial (\nabla u)= \nabla(\theta(|u|^2))$, $\theta\in C^2, \; \theta(0)=0$, then
\eqref{gNLS} has a unique classical global solution $u\in (C_w \cap
L^\infty)(\mathbb{R}, H^s) \cap C(\mathbb{R}, H^{s-1}) \cap
L^2(\mathbb{R}; H^{s-1}_{2n/(n-2)})$. The main tools used in
\cite{Oz-Zh} are  the gauge transform techniques, the energy method
together with the endpoint Strichartz estimates.

The second way consists in using the $X^{s,b}$-like spaces, see
Bourgain \cite{Bour} and it has been developed by many authors (see
\cite{Be-Ta,CKSTT,Grun} and references therein). This method depends
on both the dispersive property of the linear equation and the
structure of the nonlinearities, which is very useful for the lower
regularity initial data.

The third method is to mainly use the dispersive smooth effects of
the linear Schr\"odinger equation, see Kenig, Ponce and Vega
\cite{KePoVe1,KePoVe2}. The crucial point is that the Schr\"odinger
group has the following locally smooth effects ($n\ge 2$):
\begin{align}
& \sup_{\alpha\in \mathbb{Z}^n}\|e^{{\rm i}t \Delta}
u_0\|_{L^2_{t,x} (\mathbb{R}\times Q_\alpha)}  \lesssim
\|u_0\|_{\dot H^{-1/2}},
\label{SSM-1}\\
& \sup_{\alpha\in \mathbb{Z}^n} \left\|\nabla \int^t_0 e^{{\rm
i}(t-s) \Delta} f(s) ds \right\|_{L^2_{t,x} (\mathbb{R}\times
Q_\alpha)} \lesssim \sum_{\alpha\in \mathbb{Z}^n} \|f\|_{L^2_{t,x}
(\mathbb{R}\times Q_\alpha)}, \label{SSM-2}
\end{align}
where $Q_\alpha$ is the unit cube with center at $\alpha$. Estimate
\eqref{SSM-2} contains one order smooth effect, which can be used to
control the derivative terms in the nonlinearities. Such smooth
effect estimates are also adapted to the non-elliptic Schr\"odinger
group, i.e., \eqref{SSM-1} and \eqref{SSM-2} still hold if we
replace $e^{{\rm i}t \Delta} $ by $e^{{\rm i}t \Delta_\pm}$. Some
earlier estimates related to \eqref{SSM-1} were due to Constantin
and Saut \cite{Co-Sa}, Sj\"olin \cite{Sjol} and Vega \cite{Vega}. In
\cite{KePoVe1,KePoVe2}, the local well posedness of \eqref{gNLS} in
both elliptic and non-elliptic cases was established for
sufficiently smooth large Cauchy data ($m\ge 1$, $u_0\in H^s$ with
$s>n/2$ large enough). Moreover, they showed that the solutions are
almost global if the initial data are sufficiently small, i.e., the
maximal existing time of solutions tends to infinity as initial data
tends to 0. Recently, the local well posedness results have been
generalized to the quasi-linear (ultrahyperbolic) Schr\"odinger
equations, see \cite{KePoVe3,KePoRoVe}.  As far as the authors can
see, the existence of the scattering operators for Eq. \eqref{gNLS}
and the global well posedness of \eqref{gNLS} in the non-elliptic
cases are unknown.

\subsection{Main results}
In this paper, we mainly apply the third  method to study the global
well posedness and the existence of the scattering operators of
\eqref{gNLS} in both the elliptic and non-elliptic cases with small
data in $B^s_{2,1}$, $s> 3/2+n/2$. We now state our main results,
the notations used in this paper can be found in Sections
\ref{Notations} and \ref{functionspace}.

\begin{thm} \label{GWP-nD}
Let $n\ge 2$ and $s>n/2 +3/2$. Let $F(z)$ be as in \eqref{poly} with
$2+4/n \le m\le M<\infty $. We have the following results.

 {\rm
(i)} If $ \|u_0\|_{B^{s}_{2,1}} \le \delta$ for $n\ge 3$, and $
\|u_0\|_{B^{s}_{2,1} \cap \dot H^{-1/2}} \le \delta $ for $n=2$,
where $\delta>0$ is a suitably small number,  then \eqref{gNLS} has
a unique global solution $ u\in C(\mathbb{R}, \; B^{s}_{2,1} ) \cap
X_0, $ where
\begin{align}
X_0=\left\{u \; : \;
\begin{array}{l}
\|D^\beta u \|_{\ell^{1,s-1/2}_\triangle \ell^\infty_\alpha  (L^2_{t,x}(\mathbb{R}\times Q_\alpha))} \lesssim \delta, \ \ |\beta|\le 1 \\
\|D^\beta u \|_{\ell^{1,s-1/2}_\triangle \ell^{2+4/n}_\alpha
(L^\infty_{t,x} \cap (L^{2m}_t L^\infty_x) (\mathbb{R}\times
Q_\alpha))} \lesssim \delta, \ \  |\beta|\le 1
\end{array}
\right\}.
 \label{X0-space}
\end{align}
Moreover, for $n\ge 3$, the scattering operator of Eq. \eqref{gNLS}
carries the ball $\{u:\, \|u\|_{B^s_{2,1}} \le \delta\}$ into
$B^s_{2,1}$.

 {\rm (ii)} If $s+1/2 \in \mathbb{N}$ and $ \|u_0\|_{H^{s}}
\le \delta$ for $n\ge 3$, and $ \|u_0\|_{H^{s} \cap \dot H^{-1/2}}
\le \delta $ for $n=2$, where $\delta>0$ is a suitably small number,
then \eqref{gNLS} has a unique global solution  $ u\in C(\mathbb{R},
\; H^{s} ) \cap X, $ where
\begin{align}
X=\left\{u \; : \;
\begin{array}{l}
\|D^\beta u \|_{\ell^\infty_\alpha  (L^2_{t,x}(\mathbb{R}\times Q_\alpha))} \lesssim \delta, \;  |\beta|\le s+1/2 \\
\|D^\beta u \|_{\ell^{2+4/n}_\alpha (L^\infty_{t,x} \cap (L^{2m}_t
L^\infty_x) (\mathbb{R}\times Q_\alpha))} \lesssim \delta, \;
|\beta|\le 1
\end{array}
\right\}.
 \label{X-space}
\end{align}
Moreover, for $n\ge 3$,  the scattering operator of Eq. \eqref{gNLS}
carries the ball $\{u:\, \|u\|_{H^s} \le \delta\}$ into $H^s$.
\end{thm}
 We now illustrate the proof of
(ii) in Theorem \ref{GWP-nD}. Let us consider the equivalent
integral equation
\begin{align}
u(t) =S(t)u_0 -  {\rm i}\mathscr{A} F(u, \bar{u}, \nabla u, \nabla
\bar{u} ), \label{int-equat}
\end{align}
where
\begin{align} S(t): = e^{{\rm i} t\Delta_\pm}, \ \ \  \mathscr{A}f: =
\int^t_0 e^{{\rm i}(t-s) \Delta_\pm} f(s) ds. \label{S(t)}
\end{align}
If one applies the local smooth effect estimate \eqref{SSM-2} to
control the derivative terms in the nonlinearities, then the working
space should contains the space
$\ell^\infty_\alpha(L^2_{t,x}(\mathbb{R}\times Q_\alpha))$.  For
simplicity, we consider the case $ F(u, \bar{u}, \nabla u, \nabla
\bar{u} )$ $ = (\partial_{x_1} u)^{\nu+1}$. By \eqref{SSM-1} and
\eqref{SSM-2}, we immediately have
\begin{align}
\|\nabla u\|_{\ell^\infty_\alpha (L^2_{t,x} (\mathbb{R}\times
Q_\alpha))} & \lesssim \|u_0\|_{H^{1/2}} +\sum_{\alpha\in
\mathbb{Z}^n} \|(\partial_{x_1} u)^{\nu+1}\|_{L^2_{t,x}
(\mathbb{R}\times
Q_\alpha)} \nonumber\\
& \lesssim \|u_0\|_{H^{1/2}} + \|\nabla u\|_{\ell^\infty_\alpha
(L^2_{t,x} (\mathbb{R}\times Q_\alpha))}  \|\nabla
u\|^{\nu}_{\ell^\nu_\alpha (L^\infty_{t,x} (\mathbb{R}\times
Q_\alpha))}. \label{infty-2}
\end{align}
Hence, one needs to control $ \|\nabla u\|_{\ell^\nu_\alpha
(L^\infty_{t,x} (\mathbb{R}\times Q_\alpha))}$. In
\cite{KePoVe1,KePoVe2}, it was shown that for $\nu=2$,
\begin{align}
 \|S(t) u_0 \|_{\ell^2_\alpha (L^\infty_{t,x} ([0,T]\times Q_\alpha))} \le
C(T) \|u_0\|_{H^s}, \quad s>n/2+2. \label{Max-Hs}
\end{align}
In the elliptic case \eqref{Max-Hs} holds for $s>n/2$.
\eqref{Max-Hs} is a time-local version which prevents us to get the
global existence of solutions. So, it is natural to ask if there is
a time-global version for the estimates of the maximal function. We
can get the following
\begin{align}
 \|S(t) u_0 \|_{\ell^\nu_\alpha (L^\infty_{t,x} (\mathbb{R}\times Q_\alpha))} \le
C \|u_0\|_{H^s}, \quad s>n/2, \ \ \nu\ge 2+4/n. \label{Max-Hs1}
\end{align}
Applying \eqref{Max-Hs1}, we have for any $s>n/2$,
\begin{align}
\|\nabla u\|_{\ell^\nu_\alpha (L^\infty_{t,x} (\mathbb{R}\times
Q_\alpha))} & \lesssim \|\nabla u_0\|_{H^{s}} + \|\nabla
(\partial_{x_1}u)^{1+\nu}\|_{L^1(\mathbb{R}, H^s(\mathbb{R}^n))}.
\label{Max-L1}
\end{align}
One can get, say for $s=[n/2]+1$,
\begin{align}
 \|\nabla
(\partial_{x_1}u)^{1+\nu}\|_{L^1(\mathbb{R}, H^s(\mathbb{R}^n))} &
\lesssim  \sum_{|\beta|\le s+2}\|D^\beta u\|_{\ell^\infty_\alpha
(L^2_{t,x} (\mathbb{R}\times Q_\alpha))} \|\nabla
u\|^{\nu}_{\ell^\nu_\alpha (L^{2\nu}_t L^\infty_{x}
(\mathbb{R}\times Q_\alpha))}. \label{L1contr}
\end{align}
Hence, we need to further estimate $\|D^\beta
u\|_{\ell^\infty_\alpha (L^2_{t,x} (\mathbb{R}\times Q_\alpha))}$
for all $|\beta| \le s+2$ and $\|\nabla u\|_{\ell^\nu_\alpha
(L^{2\nu}_t L^\infty_{x} (\mathbb{R}\times Q_\alpha))}$. We can
conjecture that a similar estimate to \eqref{infty-2} holds:
\begin{align}
& \sum_{|\beta|\le s+2}\|D^\beta u\|_{\ell^\infty_\alpha (L^2_{t,x}
(\mathbb{R}\times Q_\alpha))} \nonumber\\
&  \lesssim \|u_0\|_{H^{s+ 3/2}} +   \sum_{|\beta|\le s+2}\|D^\beta
u\|_{\ell^\infty_\alpha (L^2_{t,x} (\mathbb{R}\times Q_\alpha))}
\|\nabla u\|^{\nu}_{\ell^\nu_\alpha (L^\infty_{t,x}
(\mathbb{R}\times Q_\alpha))}. \label{infty-2-1}
\end{align}
Finally, for the estimate of $\|\nabla u\|_{\ell^\nu_\alpha
(L^{2\nu}_t L^\infty_{x} (\mathbb{R}\times Q_\alpha))}$, one needs
the following
\begin{align}
 \|S(t) u_0 \|_{\ell^\nu_\alpha (L^{2\nu}_tL^\infty_{x} (\mathbb{R}\times Q_\alpha))} \le
C \|u_0\|_{H^{s-1/\nu}}, \quad s>n/2, \ \ \nu\ge 2+4/n.
\label{Max-Hs2}
\end{align}
Using \eqref{Max-Hs2},  the estimate of $\|\nabla
u\|_{\ell^\nu_\alpha (L^{2\nu}_t L^\infty_{x} (\mathbb{R}\times
Q_\alpha))}$ becomes easier than that of $\|\nabla
u\|_{\ell^\nu_\alpha (L^\infty_{t, x} (\mathbb{R}\times
Q_\alpha))}$. Hence, the solution has a self-contained behavior by
using the spaces $\ell^\infty_\alpha (L^2_{t,x} (\mathbb{R}\times
Q_\alpha)), \; \ell^\nu_\alpha (L^\infty_{t,x} (\mathbb{R}\times
Q_\alpha))$ and $\ell^\nu_\alpha (L^{2\nu}_t L^\infty_{x}
(\mathbb{R}\times Q_\alpha)) $. We will give the details of the
estimates \eqref{Max-Hs1} and \eqref{Max-Hs2} in Section
\ref{Max-Funct}. The nonlinear mapping estimates as in
\eqref{L1contr} and \eqref{infty-2-1} will be given in Section
\ref{proof-gwp-nd}.

Next, we use the frequency-uniform decomposition method developed in
\cite{Wa1,WaHe,WaHu} to consider the case of initial data in
modulation spaces $M^s_{2,1}$, which is the low regularity version
of Besov spaces $B^{n/2+s}_{2,1}$, i.e., $B^{n/2+s}_{2,1} \subset
M^s_{2,1}$ is a sharp embedding and $M^s_{2,1}$ has only $s$-order
derivative regularity (see \cite{SuTo,Toft,WaHe}, for the final
 result, see \cite{WaHu}). We have the following local well
posedness result with small rough initial data:
\begin{thm} \label{LWP-Mod}
Let $n\ge 2$. Let $F(z)$ be as in \eqref{poly} with $2 \le m\le
M<\infty $. Assume that $\|u_0\|_{M^{2}_{2,1}} \le \delta$ for $n\ge
3$, and $\|u_0\|_{M^{2}_{2,1} \cap \dot H^{-1/2}} \le \delta $ for
$n=2$, where $\delta>0$ is sufficiently small.  Then there exists a
$T:=T(\delta)>0$ such that \eqref{gNLS} has a unique local solution
$ u\in C([0,T], \; M^{2}_{2,1} ) \cap Y, $ where
\begin{align}
Y=\left\{u \; : \;
\begin{array}{l}
\|D^\beta u \|_{\ell^{1,3/2}_\Box \ell^\infty_\alpha
(L^2_{t,x}([0,T]\times Q_\alpha))} \lesssim \delta, \;  |\beta|\le 1
\\
 \|D^\beta u \|_{\ell^1_\Box \ell^{2}_\alpha (L^\infty_{t,x} ([0,T]
\times Q_\alpha))} \lesssim \delta, \; |\beta|\le 1
\end{array}
\right\}.
 \label{Y-Mod}
\end{align}
Moreover, $\lim_{\delta \searrow 0} T(\delta)=\infty$.
\end{thm}
The following is a global well posedness result with Cauchy data in
modulation spaces $M^s_{2,1}$:
\begin{thm} \label{GWP-Mod}
Let $n\ge 2$. Let $F(z)$ be as in \eqref{poly} with $2+4/n \le m\le
M<\infty $. Let $s>3/2 + (n+2)/m$.  Assume that
$\|u_0\|_{M^{s}_{2,1}} \le \delta$ for $n\ge 3$, and
$\|u_0\|_{M^{s}_{2,1} \cap \dot H^{-1/2}} \le \delta$ for $n=2$,
where  $\delta>0$ is a suitably small number.  Then  \eqref{gNLS}
has a unique global solution $ u\in C(\mathbb{R}, \; M^{s}_{2,1} )
\cap Z, $ where
\begin{align}
Z=\left\{u \; : \;
\begin{array}{l}
\|D^\beta u \|_{\ell^{1, s-1/2}_\Box \ell^\infty_\alpha
(L^2_{t,x}(\mathbb{R}
\times Q_\alpha))} \lesssim \delta, \;  |\beta|\le 1 \\
\|D^\beta u \|_{\ell^1_\Box\ell^{m}_\alpha (L^\infty_{t,x}
\cap(L^{2m}_t L^\infty_x) (\mathbb{R} \times Q_\alpha))} \lesssim
\delta, \; |\beta|\le 1
\end{array}
\right\}.
 \label{Z-Mod}
\end{align}
Moreover, for $n\ge 3$, the scattering operator of Eq. \eqref{gNLS}
carries the ball $\{u:\, \|u\|_{M^s_{2,1}} \le \delta\}$ into
$M^s_{2,1}$.
\end{thm}
Finally, we consider one spatial dimension case.  Denote
\begin{align}
s_\kappa= \frac{1}{2}- \frac{2}{\kappa}, \quad  \tilde{s}_\nu=
\frac{1}{2}- \frac{1}{\nu}. \label{Index}
\end{align}
\begin{thm} \label{GWP-1D}
Let $n=1$, $M \ge m \ge 4$,
 $u_0\in  \dot B^{1+
\tilde{s}_{M}}_{2, 1} \cap \dot B^{s_m}_{2,1}$. Assume that there
exists a small $\delta>0$ such that $ \|u_0\|_{\dot B^{1+
\tilde{s}_{M}}_{2, 1} \cap \dot B^{s_m}_{2,1}} \le \delta.$ Then
\eqref{gNLS} has a unique global solution $u\in X=\{u \in
\mathscr{S}'( \mathbb{R}^{1+1}):  \|u\|_{X} \lesssim \delta \}$,
where
\begin{align}
\|u\|_X & =\sup_{s_m\le s \le \tilde{s}_{M}} \sum_{i=0,1} \sum_{j\in
\mathbb{Z}} |\!|\!|\partial^i_x \triangle_j u |\!|\!|_s \ \ for \ \
m >4, \nonumber\\
\|u\|_X & =\sum_{i=0,1} \big(\|\partial^i_x u\|_{L^\infty_tL^2_x
\,\cap\,  L^{6}_{x,t} } + \sup_{\tilde{s}_m\le s \le \tilde{s}_{M}}
\sum_{j\in \mathbb{Z}} |\!|\!|\partial^i_x \triangle_j u |\!|\!|_s
\big) \ \ for \ \
m=4, \nonumber\\
|\!|\!|\triangle_j v & |\!|\!|_s   :=  2^{s j} (\|\triangle_j v
\|_{L^\infty_tL^2_x \,\cap\,  L^{6}_{x,t}}   + 2^{j/2}
\|\triangle_j v\|_{L^\infty_x L^2_t})\nonumber\\
& \ \ \ \ \ \  + 2^{(s-\tilde{s}_m) j}\|\triangle_j
v\|_{L_{x}^{m}L_{t}^{\infty}} +2^{(s-\tilde{s}_{M})j} \|\triangle_j
v\|_{L_x^{M}L_t^\infty}. \label{X-gNLS}
\end{align}
\end{thm}
Recall that the norm on homogeneous Besov spaces $\dot B^s_{2,1}$
can be defined in the following way:
\begin{align}
\|f\|_{\dot B^s_{2,1}} =\sum^\infty_{j=-\infty} 2^{sj}
\left(\int^{2^{j+1}}_{2^j} |\mathscr{F}f(\xi)|^2 d\xi \right)^{1/2}.
\label{hBesov}
\end{align}

\subsection{Remarks on main results}

It seems that the regularity assumptions on initial data are not
optimal in Theorems \ref{GWP-nD}--\ref{GWP-Mod}, but Theorem
\ref{GWP-1D} presents the sharp regularity condition to the initial
data. To illustrate the relation between the regularity index and
the nonlinear power, we consider a simple cases of \eqref{gNLS}:
\begin{align}
 & {\rm i} u_t + \Delta_\pm u =  u_{x_1}
^{\nu}, \ \ u(0) =\phi.  \label{dNLS}
\end{align}
Eq. \eqref{dNLS} is invariant under the scaling $u \to u_\lambda=
\lambda^{(2-\nu)/(\nu-1)} u(\lambda^2 t, \lambda x)$ and moreover,
\begin{align}
\|\phi\|_{\dot H^{s}(\mathbb{R}^n)} = \|u_\lambda (0, \cdot)\|_{\dot
H^{s}(\mathbb{R}^n)}, \quad s= 1+
\tilde{s}_{\nu-1}:=1+n/2-1/(\nu-1). \label{invariant}
\end{align}
From this point of view, we say that $s =1+ \tilde{s}_{\nu-1}$ is
the critical regularity index of \eqref{dNLS}. In \cite{MR}, Molinet
and Ribuad showed that \eqref{dNLS} is ill-posed in one spatial
dimension in the sense  if $s_1\not= \tilde{s}_{\nu-1}+ 1$,  the
flow map of equation \eqref{dNLS} $\phi \rightarrow u$ (if it
exists) is not of class $C^{\nu}$ from
$\dot{B}^{s_1}_{2,1}(\mathbb{R})$ to $C ([0, \infty),
\dot{B}^{s_1}_{2,1} (\mathbb{R}) )$ at the origin $\phi=0$. For each
term in the polynomial nonlinearity $F(u, \bar{u}, \nabla u, \nabla
\bar{u})$ as in \eqref{poly}, we easily see that the critical index
$s$ can take any critical index between $s_m$ and $1+\tilde{s}_M$.
So, our Theorem \ref{GWP-1D} give sharp result in the case $m\ge 4$.
On the other hand, Christ \cite{Chr} showed that in the case
$\nu=2$, $n=1$, for any $s\in \mathbb{R}$, there exist initial data
in $H^s$ with arbitrarily small norm, for which the solution attains
arbitrarily large norm after an arbitrarily short time (see also
\cite{MJ}). From Christ's result together with Theorems
\ref{GWP-1D}, we can expect that there exists $m_0>1$ (might be
non-integer) so that for $\nu-1\ge m_0$, $s =1+ \tilde{s}_{\nu-1}$
is the minimal regularity index to guarantee the well posedness of
\eqref{dNLS}, at least for the local solutions and small data global
solutions in $H^s$. However, it is not clear for us how to find the
exact value of $m_0$ even in one spatial dimension.

However, in higher spatial dimensions, it seems that $1/2+1/M$-order
derivative regularity is lost in Theorem \ref{GWP-nD} and we do not
know how to attain the regularity index $s\ge 1+ \tilde{s}_M$.

In two dimensional case, if $\Delta_\pm= \Delta$ and the initial
value $u_0$ is a radial function, we can remove the condition
$u_0\in \dot H^{-1/2}$, $\|u_0\|_{\dot H^{-1/2}} \le \delta$ by
using the endpoint Strichartz estimates as in the case $n\ge 3$.

Considering the nonlinearity $F(u, \nabla u)= (1-|u|^2)^{-1} |\nabla
u|^{2k} u$,  Theorem 1.2 holds for the case $k\ge 1$. Theorems
\ref{GWP-nD} and \ref{GWP-Mod} hold for the case $k\ge 2$. Since
$(1-|u|^2)^{-1}= \sum^\infty_{k=0} |u|^{2k}$,  one easily sees that
we can use the same way as in the proof of our main results to
handle this kind of nonlinearity.

\subsection{Notations} \label{Notations}

Throughout this paper,  we will always use the following notations.
$\mathscr{S}(\mathbb{R}^n)$ and $\mathscr{S}'(\mathbb{R}^n)$ stand
for the Schwartz space and its dual space, respectively. We denote
by $L^p(\mathbb{R}^n)$ the Lebesgue space, $\|\cdot\|_p:=
\|\cdot\|_{L^p(\mathbb{R}^n)}$. The Bessel potential space is
defined by $H^s_p(\mathbb{R}^n): =(I-\Delta)^{-s/2}
L^p(\mathbb{R}^n)$, $H^s(\mathbb{R}^n)=H^s_2 (\mathbb{R}^n)$, $\dot
H^s(\mathbb{R}^n) =(-\Delta)^{-s/2}
L^2(\mathbb{R}^n)$.\footnote{$\mathbb{R}^n$ will be omitted in the
definitions of various function spaces if there is no confusion.}
For any quasi-Banach space $X$, we denote by $X^*$ its dual space,
by $L^p(I, X)$ the Lebesgue-Bochner space, $\|f\|_{L^p(I,X)}:=
(\int_I \|f(t)\|^p_X dt)^{1/p}$.  If $X=L^r(\Omega)$, then we write
$L^p(I, L^r(\Omega))= L^p_tL^r_x(I\times \Omega)$ and
$L^p_{t,x}(I\times \Omega)= L^p_tL^p_x(I\times \Omega)$.   Let
$Q_\alpha$ be the unit cube with center at $\alpha\in \mathbb{Z}^n$,
i.e., $Q_\alpha=\alpha+ Q_0, Q_0= \{x=(x_1,...x_n): -1/2\le x_i<
1/2\}.$  We also needs the function spaces $\ell^q_\alpha
(L^p_tL^r_{x}(I\times Q_\alpha))$,
$$
\|f\|_{\ell^q_\alpha (L^p_{t}L^r_{x}(I\times Q_\alpha))}:=
\left(\sum_{\alpha\in \mathbb{Z}^n} \|f\|^q_{L^p_{t}L^r_{x}(I\times
Q_\alpha)} \right)^{1/q}.
$$
We denote by $\mathscr{F}$ ($\mathscr{F}^{-1}$) the (inverse)
Fourier transform for the spatial variables; by $\mathscr{F}_t$
($\mathscr{F}^{-1}_{t}$) the (inverse) Fourier transform for the
time variable and by $\mathscr{F}_{t,x}$ ($\mathscr{F}^{-1}_{t,x}$)
the (inverse) Fourier transform for both time and spatial variables,
respectively.  If there is no explanation, we always denote by
$\varphi_k(\cdot)$ the dyadic decomposition functions as in
\eqref{dyadic-funct}; and by $\sigma_k(\cdot)$ the uniform
decomposition functions as in \eqref{Mod.2}.  $u \star v$ and $u*v$
will stand for the convolution on time and on spatial variables,
respectively, i.e.,
$$
(u\star v) (t,x)=
\int_{\mathbb{R}} u(t-\tau,x) v(\tau,x)d\tau, \ \  (u* v) (t,x)=
\int_{\mathbb{R}^n} u(t,x-y) v(t,y)dy.
$$
$\mathbb{R}, \mathbb{N}$ and $ \mathbb{Z}$ will stand for the sets
of reals, positive integers and integers, respectively.  $c<1$,
$C>1$ will denote positive universal constants, which can be
different at different places. $a\lesssim b$ stands for $a\le C b$
for some constant $C>1$, $a\sim b$ means that $a\lesssim b$ and
$b\lesssim a$.  We denote by $p'$ the dual number of $p \in
[1,\infty]$, i.e., $1/p+1/p'=1$. For any $a>0$, we denote by $[a]$
the minimal integer that is larger than or equals to $a$. $B(x,R)$
will denote the ball in $\mathbb{R}^n$ with center $x$ and radial
$R$.
\subsection{Besov and modulation spaces} \label{functionspace}
Let us recall that Besov spaces $B^s_{p,q}:=B^s_{p,q}(\mathbb{R}^n)$
are defined as follows (cf. \cite{BL,Tr}). Let $\psi: \mathbb{R}^n
\to [0,1]$ be a smooth radial bump function adapted to the ball
$B(0,2)$:
\begin{align}
\psi (\xi)=\left\{
\begin{array}{ll}
1, & |\xi|\le 1,\\
{\rm smooth}, & |\xi|\in [1,2],\\
 0, & |\xi|\ge 2.
\end{array}
\right.
 \label{cutoff}
\end{align}
We write $\delta(\cdot):= \psi(\cdot)-\psi(2\,\cdot)$ and
\begin{align}
\varphi_j:= \delta(2^{-j}\cdot) \ \ {\rm for} \ \  j\ge 1; \quad
\varphi_0:=1- \sum_{j\ge 1} \varphi_j.  \label{dyadic-funct}
\end{align}

We say that $ \triangle_j := \mathscr{F}^{-1} \varphi_j \mathscr{F},
\quad j\in \mathbb{N} \cup \{0\}$ are the dyadic decomposition
operators. Beove spaces $B^s_{p,q}=B^s_{p,q}(\mathbb{R}^n)$ are
defined in the following way:
\begin{align}
B^s_{p,q} =\left \{ f\in \mathscr{S}'(\mathbb{R}^n): \;
\|f\|_{B^s_{p,q}} = \left(\sum^\infty_{j=0}2^{sjq} \|\,\triangle_j
f\|^q_p \right)^{1/q}<\infty \right\}. \label{Besov.1}
\end{align}
Now we recall the definition of modulation spaces (see
\cite{Fei2,Groh,Wa1,WaHe,WaHu}). Here we adopt an equivalent norm by
using the uniform decomposition to the frequency space.
 Let $\rho\in
\mathscr{S}(\mathbb{R}^n)$ and $\rho:\, \mathbb{R}^n\to [0,1]$ be a
smooth radial bump function adapted to the ball $B(0, \sqrt{n})$,
say $\rho(\xi)=1$ as $|\xi|\le \sqrt{n}/2$, and $\rho(\xi)=0$ as
$|\xi| \ge \sqrt{n} $. Let $\rho_k$ be a translation of $\rho$: $
\rho_k (\xi) = \rho (\xi- k), \; k\in \mathbb{Z}^n$.  We write
\begin{align}
\sigma_k (\xi)= \rho_k(\xi) \left(\sum_{k\in \mathbb{
Z}^n}\rho_k(\xi)\right)^{-1}, \quad k\in \mathbb{Z}^n. \label{Mod.2}
\end{align}
 Denote
\begin{align}
\Box_k := \mathscr{F}^{-1} \sigma_k \mathscr{F}, \quad k\in \mathbb{
Z}^n, \label{Mod.5}
\end{align}
which are said to be the frequency-uniform decomposition operators.
For any $k\in \mathbb{ Z}^n$, we write $\langle
k\rangle=\sqrt{1+|k|^2}$. Let $s\in \mathbb{R}$, $0<p,q \le \infty$.
Modulation spaces $M^s_{p,q}=M^s_{p,q}(\mathbb{R}^n)$ are defined
as:
\begin{align}
& M^s_{p,q} =\left \{ f\in \mathscr{S}'(\mathbb{R}^n): \;
\|f\|_{M^s_{p,q}} = \left(\sum_{k \in \mathbb{Z}^n} \langle
k\rangle^{sq} \|\,\Box_k f\|_p^q \right)^{1/q}<\infty \right\}.
\label{Mod.6}
\end{align}
We will use the function space $\ell^{1,s}_\Box \ell^q_\alpha
(L^p_{t}L^r_x (I\times Q_\alpha))$ which contains all of the
functions $f(t,x)$ so that the following norm is finite:
\begin{align}
\|f\|_{\ell^{1,s}_\Box \ell^q_\alpha (L^p_{t} L^r_x(I\times
Q_\alpha))}:= \sum_{k\in \mathbb{Z}^n} \langle k\rangle ^s
\left(\sum_{\alpha\in \mathbb{Z}^n} \|\Box_k f\|^q_{L^p_{t}L^r_x
(I\times Q_\alpha)} \right)^{1/q}. \label{Mod.7}
\end{align}
  Similarly, we can define the space $\ell^{1,s}_\triangle \ell^q_\alpha (L^p_{t,x}(I\times
Q_\alpha))$ with the following norm:
\begin{align}
\|f\|_{\ell^{1,s}_\triangle \ell^q_\alpha (L^p_{t} L^r_x (I\times
Q_\alpha))}:= \sum^\infty_{j=0} 2^{sj} \left(\sum_{\alpha\in
\mathbb{Z}^n} \|\triangle_j f\|^q_{L^p_{t}L^r_x (I\times Q_\alpha)}
\right)^{1/q}. \label{Besov.7}
\end{align}
A special case is $s=0$, we write  $\ell^{1,0}_\Box \ell^q_\alpha
(L^p_{t}L^r_x (I\times Q_\alpha))= \ell^{1}_\Box \ell^q_\alpha
(L^p_{t}L^r_x (I\times Q_\alpha))$ and $\ell^{1,0}_\triangle
\ell^q_\alpha (L^p_{t}L^r_x(I\times Q_\alpha)) = \ell^{1}_\triangle
\ell^q_\alpha (L^p_{t}L^r_x (I\times Q_\alpha))$.

 The rest of this paper is organized as follows. In
Section \ref{Max-Funct} we give the details of the estimates for the
maximal function in certain function spaces. Section
\ref{global-local-estimate} is devoted to considering the spatial
local versions for the Strichartz estimates and giving some remarks
on the estimates of the local smooth effects. In Sections
\ref{proof-gwp-nd}--\ref{proof-gwp-1D} we prove our main Theorems
\ref{GWP-nD}--\ref{GWP-1D}, respectively.

\section{Estimates for the maximal function} \label{Max-Funct}
\subsection{Time-local version}

Recall that $S(t)= e^{- {\rm i}t \triangle_\pm }= \mathscr{F}^{-1}
e^{ {\rm i}t |\xi|^2_\pm } \mathscr{F}$, where
\begin{align}
|\xi|^2_\pm  = \sum^n_{j=1} \varepsilon_j \xi^2_j, \quad
\varepsilon_j = \pm 1. \label{symbol}
\end{align}
Kenig, Ponce and Vega \cite{KePoVe1} showed the following maximal
function estimate:
\begin{align}
\left( \sum_{\alpha\in \mathbb{Z}^n}
\|S(t)u_0\|^2_{L^\infty_{t,x}([0,T]\times Q_\alpha)} \right)^{1/2}
\lesssim C(T) \|u_0\|_{H^s}, \label{max-0}
\end{align}
where $s\ge 2+n/2$. If $S(t)=  e^{- {\rm i}t \triangle }$, then
\eqref{max-0} holds for $s>n/2$, $C(T)=(1+T)^s$.  Using the
frequency-uniform decomposition method, we can get the following

\begin{prop} \label{prop-max-1}
There exists a constant $C(T)>1$ which depends only on $T$ and $n$
such that
\begin{align}
\sum_{k\in \mathbb{Z}^n}\left( \sum_{\alpha\in \mathbb{Z}^n}
\|\,\Box_k S(t)u_0\|^2_{L^\infty_{t,x}([0,T]\times Q_\alpha)}
\right)^{1/2} \le C(T)  \|u_0\|_{M^{1/2}_{2,1}}, \label{max-1}
\end{align}
\end{prop}
In particular, for any $s>(n+1)/2$,
\begin{align}
\left( \sum_{\alpha\in \mathbb{Z}^n} \|\,
S(t)u_0\|^2_{L^\infty_{t,x}([0,T]\times Q_\alpha)} \right)^{1/2} \le
C(T)  \|u_0\|_{H^{s}}. \label{max-1-note}
\end{align}

\noindent {\bf Proof.} By the duality, it suffices to prove that
\begin{align}
\int^T_0 (S(t)u_0, \psi(t) ) dt  \lesssim \|u_0\|_{M^{1/2}_{2,1}}
\sup_{k\in \mathbb{Z}^n}\left( \sum_{\alpha\in \mathbb{Z}^n}
\|\Box_k \psi(t)\|_{L^1_{t,x}([0,T]\times Q_\alpha)}^2
\right)^{1/2}. \label{max-2}
\end{align}
Since $(M^{1/2}_{2,1})^*= M^{-1/2}_{2,\infty}$, we have
\begin{align}
\int^T_0 (S(t)u_0, \psi(t) ) dt  \le \|u_0\|_{M^{1/2}_{2,1}} \left\|
\int^T_0 S(-t) \psi(t) dt \right\|_{M^{-1/2}_{2,\infty}}.
\label{max-3}
\end{align}
Recalling that $\|f\|_{M^{-1/2}_{2,\infty}} =\sup_{k\in
\mathbb{Z}^n} \langle k\rangle^{-1/2} \|\Box_k f\|_2$, we need to
estimate
\begin{align}
& \left\| \Box_k \int^T_0 S(-t) \psi(t) dt \right\|^2_2
\nonumber\\
& = \int^T_0 \left( \Box_k \psi(t), \;  \int^T_0 S(t-\tau) \Box_k
\psi(\tau) d\tau \right) dt \nonumber\\
& \le \sum_{\alpha \in \mathbb{Z}^n} \|\Box_k
\psi\|_{L^1_{t,x}([0,T]\times Q_\alpha)}  \left\|\int^T_0 S(t-\tau)
\Box_k
\psi(\tau) d\tau \right\|_{L^\infty_{t,x}([0,T]\times Q_\alpha)} \nonumber\\
& \le  \|\Box_k \psi\|_{\ell^2_\alpha (L^1_{t,x}([0,T]\times
Q_\alpha))}  \left\|\int^T_0 S(t-\tau) \Box_k \psi(\tau) d\tau
\right\|_{\ell^2_\alpha (L^\infty_{t,x}([0,T]\times Q_\alpha))} .
\label{max-4}
\end{align}
If one can show that
\begin{align}
  \left\|\int^T_0 S(t-\tau) \Box_k \psi(\tau) d\tau
\right\|_{\ell^2_\alpha (L^\infty_{t,x}([0,T]\times Q_\alpha))}
\lesssim C(T)\langle k\rangle  \|\Box_k \psi\|_{\ell^2_\alpha
(L^1_{t,x}([0,T]\times Q_\alpha))}, \label{max-5}
\end{align}
then from \eqref{max-3}--\eqref{max-5} we obtain that \eqref{max-2}
holds. Denote
\begin{align}
\Lambda:= \{\ell \in \mathbb{Z}^n: {\rm supp} \,\sigma_\ell \cap
{\rm supp} \,\sigma_0 \not= \varnothing \}.  \label{orthogonal}
\end{align}
In  the following we show \eqref{max-5}. In view of Young's
inequality, we have
\begin{align}
&  \left\|\int^T_0 S(t-\tau) \Box_k \psi(\tau) d\tau
\right\|_{L^\infty_{t,x}([0,T]\times Q_\alpha)} \nonumber\\
& \lesssim  \sum_{\ell \in \Lambda} \left\|\int^T_0 S(t-\tau)
\Box_{k+\ell} \Box_k \psi(\tau) d\tau
\right\|_{L^\infty_{t,x}([0,T]\times Q_\alpha)} \nonumber\\
& = \sum_{\ell \in \Lambda} \left\|\int^T_0  [\mathscr{F}^{-1} (
e^{{\rm i} (t-\tau) |\xi|^2_\pm} \sigma_{k+\ell})]*
 \Box_k \psi(\tau) d\tau
\right\|_{L^\infty_{t,x}([0,T]\times Q_\alpha)} \nonumber\\
& \le \sum_{\ell \in \Lambda} \sum_{\beta\in \mathbb{Z}^n}
\left\|\mathscr{F}^{-1} ( e^{{\rm i} t |\xi|^2_\pm} \sigma_{k+\ell})
\right\|_{L^\infty_{t,x}([-T,T]\times Q_\beta)} \| \Box_k
\psi\|_{L^1_{t,x}([0,T]\times (Q_\alpha-Q_\beta))}.  \label{max-6}
\end{align}
From \eqref{max-6} and Minkowski's inequality that
\begin{align}
&  \left\|\int^T_0 S(t-\tau) \Box_k \psi(\tau) d\tau
\right\|_{\ell^2_\alpha (L^\infty_{t,x}([0,T]\times Q_\alpha))} \nonumber\\
& \le  \sum_{\ell \in \Lambda} \sum_{\beta\in \mathbb{Z}^n}
\left\|\mathscr{F}^{-1} ( e^{{\rm i} t |\xi|^2_\pm} \sigma_{k+\ell})
\right\|_{L^\infty_{t,x}([-T,T]\times Q_\beta)} \| \Box_k
\psi\|_{\ell^2_\alpha (L^1_{t,x}([0,T]\times (Q_\alpha-Q_\beta)))}.
\label{max-7}
\end{align}
It is easy to see that
\begin{align}
 \| \Box_k
\psi\|_{\ell^2_\alpha (L^1_{t,x}([0,T]\times (Q_\alpha-Q_\beta)))}
\lesssim   \| \Box_k \psi\|_{\ell^2_\alpha (L^1_{t,x}([0,T]\times
Q_\alpha))}. \label{max-8}
\end{align}
Hence, in order to prove \eqref{max-5}, it suffices to prove that
\begin{align}
 \sum_{\beta\in \mathbb{Z}^n}
\left\|\mathscr{F}^{-1} ( e^{{\rm i} t |\xi|^2_\pm} \sigma_{k})
\right\|_{L^\infty_{t,x}([0,T]\times Q_\beta)}  \lesssim C(T)
\langle k \rangle. \label{max-9}
\end{align}
In fact, observing the following identity,
\begin{align}
 |\mathscr{F}^{-1} ( e^{{\rm i} t |\xi|^2_\pm} \sigma_{k}) |=
 |\mathscr{F}^{-1} ( e^{{\rm i} t |\xi|^2_\pm} \sigma_{0}) (\cdot + 2tk_\pm) |,  \label{max-10}
\end{align}
where $k_\pm= (\varepsilon_1 k_1,..., \varepsilon_n k_n)$, we have
\begin{align}
 \|\mathscr{F}^{-1} ( e^{{\rm i} t |\xi|^2_\pm} \sigma_{k}) \|_{L^\infty_{t,x}([0,T]\times Q_\beta)}
 \le
 \|\mathscr{F}^{-1} ( e^{{\rm i} t |\xi|^2_\pm} \sigma_{0})\|_{L^\infty_{t,x}([0,T]\times Q^*_{\beta,k})},  \label{max-11}
\end{align}
where
\begin{align*}
  Q^*_{\beta,k}=  \{x: \; x\in 2tk_\pm + Q_\beta \;\; \mbox{for some }\; t\in [0,T]\}.
\end{align*}
Denote $\Lambda_{\beta,k} = \{\beta': \; Q_{\beta'} \cap
Q^*_{\beta,k} \neq \varnothing\}$. It follows from \eqref{max-11}
that
\begin{align}
\sum_{\beta \in \mathbb{Z}^n} \|\mathscr{F}^{-1} ( e^{{\rm i} t
|\xi|^2_\pm} \sigma_{k}) \|_{L^\infty_{t,x}([0,T]\times Q_\beta)}
 \le \sum_{\beta \in \mathbb{Z}^n} \sum_{\beta' \in \Lambda_{\beta,k}}
 \|\mathscr{F}^{-1} ( e^{{\rm i} t |\xi|^2_\pm} \sigma_{0})\|_{L^\infty_{t,x}([0,T]\times Q_{\beta'})}.  \label{max-13}
\end{align}
Since each $E_{\beta,k}$ overlaps at most $O(T\langle k\rangle)$
many $Q_{\beta'}$, $\beta'\in \mathbb{Z}^n$, one can easily verify
that in the sums of the right hand side of \eqref{max-13}, each
$\|\mathscr{F}^{-1} ( e^{{\rm i} t |\xi|^2_\pm}
\sigma_{0})\|_{L^\infty_{t,x}([0,T]\times Q_{\beta'})}$ repeats at
most  $O(T \langle k\rangle)$ times. Hence, we have
\begin{align}
\sum_{\beta \in \mathbb{Z}^n} \|\mathscr{F}^{-1} ( e^{{\rm i} t
|\xi|^2_\pm} \sigma_{k}) \|_{L^\infty_{t,x}([0,T]\times Q_\beta)}
 \lesssim \langle k\rangle  \sum_{\beta \in \mathbb{Z}^n}
 \|\mathscr{F}^{-1} ( e^{{\rm i} t |\xi|^2_\pm} \sigma_{0})\|_{L^\infty_{t,x}([0,T]\times Q_{\beta})}.  \label{max-14}
\end{align}
Finally, it suffices to show that
\begin{align}
\sum_{\beta \in \mathbb{Z}^n} \|\mathscr{F}^{-1} ( e^{{\rm i} t
|\xi|^2_\pm} \sigma_{0})\|_{L^\infty_{t,x}([0,T]\times Q_{\beta})}
\le C(T).  \label{max-15}
\end{align}
Denote $\nabla_{t,x} =(\partial_t, \partial_{x_1},...,
\partial_{x_n})$. By the Sobolev  inequality,
\begin{align}
\sum_{\beta \in \mathbb{Z}^n} \|\mathscr{F}^{-1} ( e^{{\rm i} t
|\xi|^2_\pm} \sigma_{0})\|_{L^\infty_{t,x}([0,T]\times Q_{\beta})}
\lesssim & \; \sum_{\beta \in \mathbb{Z}^n} \|\mathscr{F}^{-1} (
e^{{\rm i} t |\xi|^2_\pm}
\sigma_{0})\|_{L^{2n}_{t,x}([0,T]\times Q_{\beta})}\nonumber\\
& + \sum_{\beta \in \mathbb{Z}^n} \|\nabla_{t,x} \mathscr{F}^{-1} (
e^{{\rm i} t |\xi|^2_\pm} \sigma_{0})\|_{L^{2n}_{t,x}([0,T]\times
Q_{\beta})} \nonumber\\
 = & I+II. \label{max-16}
\end{align}
By H\"older's inequality, we have
\begin{align}
II  \lesssim &   \left(\sum_{\beta \in \mathbb{Z}^n}
\left\|(1+|x|^2)^n\nabla_{t,x} \mathscr{F}^{-1} ( e^{{\rm i} t
|\xi|^2_\pm} \sigma_{0}) \right\|^{2n}_{L^{2n}_{t,x}([0,T]\times
Q_{\beta})}\right)^{1/2n} \nonumber\\
\lesssim &   \sum^n_{i=1} \left\|\mathscr{F}^{-1} (I-\Delta)^n (
e^{{\rm i} t |\xi|^2_\pm} \xi_i \sigma_{0})
\right\|_{L^{2n}_{t,x}([0,T]\times \mathbb{R}^n)}\nonumber\\
 &
+ \left\|\mathscr{F}^{-1} (I-\Delta)^n ( e^{{\rm i} t |\xi|^2_\pm}
|\xi|^2_\pm \sigma_{0}) \right\|_{L^{2n}_{t,x}([0,T]\times
\mathbb{R}^n)}\nonumber\\
\lesssim & \; C(T). \label{max-17}
\end{align}
One easily sees that $I$ has the same bound as that of $II$. The
proof of \eqref{max-1} is finished.  Noticing that $H^s\subset
M^{1/2}_{2,1}$ if $s>(n+1)/2$ (cf. \cite{Toft,WaHe,WaHu}), we
immediately have \eqref{max-1-note}. $ \hfill\Box$

\subsection{Time-global version}

Recall that we have the following equivalent norm on Besov spaces
(\cite{BL,Tr}):

\begin{lem} \label{Besovnorm}
Let $1\le p,q \le \infty$, $\sigma>0$, $\sigma\not\in \mathbb{N}$.
Then we have
\begin{align}
\|f\|_{B^\sigma_{p,q}} \sim \sum_{|\beta|\le [\sigma]} \|D^\beta
f\|_{L^p(\mathbb{R}^n)}+ \sum_{|\beta|\le [\sigma]} \left(
\int_{\mathbb{R}^n} |h|^{-n-q\{\sigma\}} \| \vartriangle_h D^\beta f
\|^q_{L^p(\mathbb{R}^n)} dh \right)^{1/q}, \label{Besov-norm-1}
\end{align}
where $\vartriangle_h f= f(\cdot+h)-f(\cdot)$, $[\sigma]$ denotes
the minimal integer that is larger than or equals to $\sigma$,
$\{\sigma\}= \sigma-[\sigma]$.
\end{lem}
Taking $p=q$ in Lemma \ref{Besovnorm}, one has that
\begin{align}
\|f\|^p_{B^\sigma_{p,p}} \sim \sum_{|\beta|\le [\sigma]} \|D^\beta
f\|^p_{L^p(\mathbb{R}^n)}+ \sum_{|\beta|\le [\sigma]}
\int_{\mathbb{R}^n} \int_{\mathbb{R}^n} \frac{| \vartriangle_h
D^\beta f (x)|^p}{|h|^{n+p\{\sigma\}}}  dxdh. \label{Besov-norm-2}
\end{align}

\begin{lem} \label{Max-Besov-control}
Let $1< p< \infty$,  $s>1/p$. Then we have
\begin{align}
\left( \sum_{\alpha\in \mathbb{Z}^n} \|u\|^p_{L^\infty_{t,x}
(\mathbb{R} \times Q_\alpha)} \right)^{1/p} \lesssim
\|(I-\partial^2_t)^{s/2} u\|_{L^p(\mathbb{R},
B^{ns}_{p,p}(\mathbb{R}^n))}. \label{Max-Besov-1}
\end{align}
\end{lem}
{\bf Proof.} We divide the proof into the following two cases.

{\it Case} 1. $ns\not\in \mathbb{N}.$  Due to $H^s_p(\mathbb{R})
\subset L^\infty(\mathbb{R})$, we have
\begin{align}
\|u\|_{L^\infty_{t,x} (\mathbb{R} \times Q_\alpha)} & \lesssim
\|(I-\partial^2_t)^{s/2} u\|_{L^\infty_x L^p_t(Q_\alpha\times
\mathbb{R})} \nonumber\\
& \le \|(I-\partial^2_t)^{s/2} u\|_{L^p_t L^\infty_x (\mathbb{R}
\times Q_\alpha )} . \label{Max-Besov-2}
\end{align}
Recalling that $\sigma_\alpha (x) \gtrsim 1$ for all $x\in Q_\alpha$
and $\alpha\in \mathbb{Z}^n$, we have from \eqref{Max-Besov-2} that
\begin{align}
\|u\|_{L^\infty_{t,x} (\mathbb{R} \times Q_\alpha)}  \le
\|(I-\partial^2_t)^{s/2} \sigma_\alpha u\|_{L^p_t L^\infty_x
(\mathbb{R} \times \mathbb{R}^n )} . \label{Max-Besov-3}
\end{align}
Since $B^{ns}_{p,p} ( \mathbb{R}^n) \subset L^\infty (
\mathbb{R}^n)$, in view of \eqref{Max-Besov-3}, one has that
\begin{align}
\|u\|_{L^\infty_{t,x} (\mathbb{R} \times Q_\alpha)}  \le
\|(I-\partial^2_t)^{s/2} \sigma_\alpha u\|_{L^p(\mathbb{R},
B^{ns}_{p,p}(\mathbb{R}^n))}. \label{Max-Besov-4}
\end{align}
For simplicity, we denote $v=(I-\partial^2_t)^{s/2}  u$.   By
\eqref{Besov-norm-2}  and \eqref{Max-Besov-4} we have
\begin{align}
\sum_{\alpha\in \mathbb{Z}^n} \|u\|^p_{L^\infty_{t,x}
(\mathbb{R}\times Q_\alpha)} \lesssim & \sum_{|\beta|\le [ns]}
\sum_{\alpha\in \mathbb{Z}^n} \int_{\mathbb{R}} \|D^\beta
(\sigma_\alpha
v)(t)\|^p_{L^p(\tilde{Q}_\alpha)} dt \nonumber\\
 & + \sum_{|\beta|\le [ns]} \sum_{\alpha\in \mathbb{Z}^n} \int_{\mathbb{R}}
\int_{\mathbb{R}^n} \int_{\mathbb{R}^n} \frac{| \vartriangle_h
D^\beta (\sigma_\alpha v)(t,x)|^p}{|h|^{n+p\{ns\}}}  dxdh dt
\nonumber\\
:= & I+II. \label{Max-Besov-5}
\end{align}
We now estimate $II$. It is easy to see that
\begin{align}
| \vartriangle_h \! D^\beta  (\sigma_\alpha v)| & \lesssim
\sum_{\beta_1+
\beta_2=\beta}  |\vartriangle_h \! (D^{\beta_1} \sigma_\alpha  D^{\beta_2} v) |  \nonumber\\
 & \leqslant \sum_{\beta_1+
\beta_2=\beta}  (|D^{\beta_1} \sigma_\alpha  (\cdot+h) \,
\vartriangle_h\! D^{\beta_2} v | + |(\vartriangle_h\! D^{\beta_1}
\sigma_\alpha) D^{\beta_2} v |). \label{Max-Besov-6}
\end{align}
Since ${\rm supp}\,\sigma_\alpha$ overlaps at most finitely many
${\rm supp}\, \sigma_\beta$ and $\sigma_\beta =\sigma_0
(\cdot-\beta)$, $\beta\in \mathbb{Z}^n$,  it follows from
\eqref{Max-Besov-6}, $|D^{\beta_1} \sigma_\alpha| \lesssim 1$ and
H\"older's inequality that
\begin{align}
II  \lesssim &   \sum_{|\beta_1|, |\beta_2|\le [ns]}
\int_{\mathbb{R}} \int_{\mathbb{R}^n} \int_{\mathbb{R}^n}
\sum_{\alpha\in \mathbb{Z}^n} |D^{\beta_1} \sigma_{\alpha}(x+h)|
\frac{| \vartriangle_h\! D^{\beta_2} v (t,x)|^p}{|h|^{n+p\{ns\}}}
dxdh dt
\nonumber\\
& +  \sum_{|\beta|\le [ns]}\sum_{\beta_1+\beta_2=\beta}
\int_{\mathbb{R}^n}    \frac{\| \vartriangle_h\! D^{\beta_1}
\sigma_0  \|^p_{L^\infty(\mathbb{R}^n)}}{|h|^{n+p\{ns\}}} dh
\nonumber\\
& \quad \times \sup_{h} \sum_{\alpha\in \mathbb{Z}^n}
\int_{\mathbb{R}} \int_{B(0,\sqrt{n}) \cup B(-h, \sqrt{n}))}
|D^{\beta_2} v(t,x+\alpha)|^p dx dt
\nonumber\\
\lesssim &   \sum_{|\beta|\le [ns]}  \int_{\mathbb{R}}
\int_{\mathbb{R}^n} \int_{\mathbb{R}^n} \frac{| \vartriangle_h\!
D^\beta v (t,x)|^p}{|h|^{n+p\{ns\}}}  dxdh dt
\nonumber\\
& +  \|\sigma_\alpha\|^p_{B^{ns}_{\infty, p}} \sum_{|\beta|\le [ns]}
\int_{\mathbb{R}} \int_{\mathbb{R}^n} |D^{\beta} v(x)|^p dx dt
\nonumber\\
\lesssim & \; \|v\|^p_{L^p(\mathbb{R}, B^{ns}_{p,p}(\mathbb{R}^n))}.
\label{Max-Besov-7}
\end{align}
Clearly, one has that
\begin{align}
I \lesssim  \|v\|^p_{L^p(\mathbb{R}, B^{ns}_{p,p}(\mathbb{R}^n))}.
\label{Max-Besov-8}
\end{align}
Collecting \eqref{Max-Besov-5}, \eqref{Max-Besov-7} and
\eqref{Max-Besov-8}, we have \eqref{Max-Besov-1}.

{\it Case} 2. $ns \in \mathbb{N}.$ One can take an $s_1< s$ such
that $s_1>1/p$ and $ns_1 \not\in \mathbb{N}$. Applying the
conclusion as in Case 1, we get the result, as desired. $\hfill\Box$

\medskip

For the semi-group $S(t)$, we have the following Strichartz estimate
(cf. \cite{Ke-Ta}):
\begin{prop}\label{Strichartz-Leb}
Let $n\ge 2$.  $2\le p, \rho \le 2n/(n-2)$ $(2\le p, \rho<\infty$ if
$n=2)$, $2/\gamma(\cdot)= n(1/2-1/\cdot)$. We have
\begin{align}
\|S(t)u_0\|_{L^{\gamma(p)}(\mathbb{R}, L^p(\mathbb{R}^n))} &
\lesssim
\| u_0\|_{L^2(\mathbb{R}^n)}, \label{Str-Leb}\\
\|\mathscr{A} F\|_{L^{\gamma(p)}(\mathbb{R}, L^p(\mathbb{R}^n))} &
\lesssim  \| F\|_{L^{\gamma(\rho)'}(\mathbb{R},
L^{\rho'}(\mathbb{R}^n))}.  \label{Str-Leb-Int}
\end{align}
\end{prop}
If $p$ and $\rho$ equal to $2n/(n-2)$, then \eqref{Str-Leb} and
\eqref{Str-Leb-Int} are said to be the endpoint Strichartz
estimates.  Using Proposition \ref{Strichartz-Leb}, we have
\begin{prop} \label{prop-max-2}
Let $p\ge 2+ 4/n:= 2^*$. For any $s>n/2$,  we have
\begin{align}
\left( \sum_{\alpha\in \mathbb{Z}^n} \|\,
S(t)u_0\|^p_{L^\infty_{t,x}(\mathbb{R} \times Q_\alpha)}
\right)^{1/p} \lesssim  \|u_0\|_{H^{s}}. \label{m-1}
\end{align}
\end{prop}
{\bf Proof.} For short, we write $\langle \partial_t
\rangle=(I-\partial^2_t)^{1/2}$.  By Lemma \ref{Max-Besov-control},
for any $s_0>1/2^*$,
\begin{align}
\left( \sum_{\alpha\in \mathbb{Z}^n} \|S(t)u_0\|^p_{L^\infty_{t,x}
(\mathbb{R} \times Q_\alpha)} \right)^{1/p} & \lesssim \left(
\sum_{\alpha\in \mathbb{Z}^n} \|S(t)u_0\|^{2^*}_{L^\infty_{t,x}
(\mathbb{R}
\times Q_\alpha)} \right)^{1/{2^*}} \nonumber\\
& \lesssim \|\langle \partial_t \rangle^{s_0}
S(t)u_0\|_{L^{2^*}(\mathbb{R}, B^{ns_0}_{2^*, 2^*}(\mathbb{R}^n))}.
\label{MB-1}
\end{align}
We have
\begin{align}
\|\langle \partial_t \rangle^{s_0}
S(t)u_0\|^{2^*}_{L^{2^*}(\mathbb{R}, B^{ns_0}_{2^*,
2^*}(\mathbb{R}^n))}= \sum_{k=0}^\infty 2^{ns_0k 2^*} \|\langle
\partial_t \rangle^{s_0} \triangle_k S(t)u_0\|^{2^*}_{L^{2^*}_{t,x}
(\mathbb{R}^{1+n})}. \label{MB-2}
\end{align}
Using the dyadic decomposition to the time-frequency, we obtain that
\begin{align}
\|\langle \partial_t \rangle^{s_0} \triangle_k
S(t)u_0\|_{L^{2^*}_{t,x}} \lesssim \sum^\infty_{j=0}  \|
\mathscr{F}^{-1}_{t,x}  \langle \tau \rangle^{s_0} \varphi_j (\tau)
\mathscr{F}_t e^{{\rm i}t|\xi|^2_\pm} \varphi_k (\xi) \mathscr{F}_x
u_0\|_{L^{2^*}_{t,x}}. \label{MB-3}
\end{align}
Noticing the fact that
\begin{align}
 (\mathscr{F}^{-1}_{t}  \langle \tau
\rangle^{s_0} \varphi_j (\tau)) \star e^{{\rm i}t|\xi|^2_\pm} = c\;
e^{{\rm i}t|\xi|^2_\pm} \varphi_j (|\xi|^2_\pm)\langle |\xi|^2_\pm
\rangle^{s_0}, \label{MB-4}
\end{align}
and using the Strichartz inequality and Plancherel's identity, one
has that
\begin{align}
\|\langle \partial_t \rangle^{s_0} \triangle_k
S(t)u_0\|_{L^{2^*}_{t,x}} & \lesssim \sum^\infty_{j=0}  \| S(t)
\mathscr{F}^{-1}_{x} \langle |\xi|^2_\pm \rangle^{s_0} \varphi_j
(|\xi|^2_\pm)  \varphi_k (\xi) \mathscr{F}_x u_0\|_{L^{2^*}_{t,x}} \nonumber\\
& \lesssim \sum^\infty_{j=0}  \| \mathscr{F}^{-1}_{x} \langle
|\xi|^2_\pm \rangle^{s_0} \varphi_j
(|\xi|^2_\pm)  \varphi_k (\xi) \mathscr{F}_x u_0\|_{L^{2}_{x}(\mathbb{R}^n)} \nonumber\\
& \lesssim  2^{2 s_0 k} \sum^\infty_{j=0}  \| \mathscr{F}^{-1}_{x}
 \varphi_j
(|\xi|^2_\pm)  \varphi_k (\xi) \mathscr{F}_x
u_0\|_{L^{2}_{x}(\mathbb{R}^n)}. \label{MB-5}
\end{align}
Combining \eqref{MB-2} and  \eqref{MB-5},  together with Minkowski's
inequality, we have
\begin{align}
& \|\langle \partial_t \rangle^{s_0} S(t)u_0\|_{L^{2^*}(\mathbb{R},
B^{ns_0}_{2^*,
2^*}(\mathbb{R}^n))} \nonumber\\
& \lesssim \sum^\infty_{j=0} \left( \sum^\infty_{k=0} 2^{(n+2)s_0k
2^*}\| \mathscr{F}^{-1}
 \varphi_j
(|\xi|^2_\pm)  \varphi_k \mathscr{F}
u_0\|^{2^*}_{L^{2}_{x}(\mathbb{R}^n)} \right)^{1/2^*} \nonumber\\
& \lesssim \sum^\infty_{j=0}  \| \mathscr{F}^{-1}
 \varphi_j
(|\xi|^2_\pm)   \mathscr{F} u_0\|_{B^{(n+2)s_0}_{2, 2^*}}.
 \label{MB-6}
\end{align}
In view of $H^{(n+2)s_0}\subset B^{(n+2)s_0}_{2,2^*}$ and H\"older's
inequality, we have for any $\varepsilon>0$,
\begin{align}
\sum^\infty_{j=0}  \| \mathscr{F}^{-1}
 \varphi_j
(|\xi|^2_\pm)   \mathscr{F} u_0\|_{B^{(n+2)s_0}_{2, 2^*}} & \lesssim
\sum^\infty_{j=0}  \| \mathscr{F}^{-1}
 \varphi_j
(|\xi|^2_\pm)   \mathscr{F} u_0\|_{H^{(n+2)s_0}} \nonumber\\
 & \lesssim
\left(\sum^\infty_{j=0}  2^{2j\varepsilon} \| \mathscr{F}^{-1}
 \varphi_j
(|\xi|^2_\pm)   \mathscr{F} u_0\|^2_{H^{(n+2)s_0}} \right)^{1/2}.
 \label{MB-7}
\end{align}
By Plancherel's identity, and ${\rm supp} \varphi_j (|\xi|^2_\pm)
\subset \{\xi:\; ||\xi|^2_\pm|\in [2^{j-1}, 2^{j+1}]\}$, we easily
see that
\begin{align}
\left(\sum^\infty_{j=0}  2^{2j\varepsilon} \| \mathscr{F}^{-1}
 \varphi_j
(|\xi|^2_\pm)   \mathscr{F} u_0\|^2_{H^{(n+2)s_0}} \right)^{1/2}
 &
\lesssim \left(\sum^\infty_{j=0}   \| \langle
|\xi|^2_\pm\rangle^\varepsilon
 \varphi_j
(|\xi|^2_\pm)   \mathscr{F} u_0\|^2_{H^{(n+2)s_0}} \right)^{1/2}
\nonumber\\
&  \lesssim \left(\sum^\infty_{j=0}   \|
 \varphi_j
(|\xi|^2_\pm)   \mathscr{F} u_0\|^2_{H^{(n+2)s_0+2\varepsilon}}
\right)^{1/2}
\nonumber\\
&  \lesssim   \| u_0\|_{H^{(n+2)s_0+2\varepsilon}} .
 \label{MB-8}
\end{align}
Taking $s_0$ such that $(n+2)s_0 +2 \varepsilon< s$, from
\eqref{MB-6}--\eqref{MB-8} we have the result, as desired.
$\hfill\Box$

\medskip

Next, we consider the estimates for the maximal function based on
the frequency-uniform decomposition method. This issue has some
relations with the Strichartz estimates in modulation spaces.
Recently, the Strichartz estimates have been generalized to various
function spaces, for instance, in the Wiener amalgam spaces
\cite{Co-Ni1,Co-Ni2}.  Recall that in \cite{WaHe}, we obtained the
following Strichartz estimate for a class of dispersive semi-groups
in modulation spaces:
\begin{align}
U(t)= \mathscr{F}^{-1} e^{{\rm i}t P(\xi)} \mathscr{F},
\label{Str-0}
\end{align}
$P(\cdot): \mathbb{R}^n\to \mathbb{R}$ is a real valued function,
which satisfies the following decay estimate
\begin{align}
\|U(t)f\|_{M^\alpha_{ p,q}} \lesssim (1+|t|)^{-\delta} \| f\|_{M_{
p',q}}, \label{Str-1}
\end{align}
where $2\le p<\infty$, $\alpha=\alpha(p)\in \mathbb{R}$,
$\delta=\delta(p)>0$, $\alpha, \delta$ are independent of $t\in
\mathbb{R}$.

\begin{prop}\label{Strichartz-Mod}
Let $U(t)$ satisfy \eqref{Str-1} and \eqref{Str-2}.  We have for any
$\gamma\ge 2 \vee (2/\delta)$,
\begin{align}
\|U(t)f\|_{L^\gamma(\mathbb{R}, M^{\alpha/2}_{p,1})} \lesssim \|
f\|_{M_{2,1}}. \label{Str-2}
\end{align}
\end{prop}
Recall that the hyperbolic Schr\"odinger semi-group $S(t)=e^{{\rm i}
t\Delta_\pm}$ has the same decay estimate as that of the elliptic
Schr\"odinger semi-group $e^{{\rm i} t\Delta}$:
\begin{align*}
\|S(t) u_0 \|_{L^\infty(\mathbb{R}^n)} \lesssim  |t|^{-n/2} \|
u_0\|_{L^1(\mathbb{R}^n))}.
\end{align*}
It follows that
\begin{align}
\|S(t)u_0\|_{M_{\infty,1}} \lesssim  |t|^{-n/2}  \| u_0\|_{M_{1,1}}.
\label{Str-3}
\end{align}
On the other hand, by Hausdorff-Young's and H\"older's inequalities
we easily calculate that
\begin{align}
\|\,\Box_k S(t)u_0\|_{L^\infty(\mathbb{R}^n)} & \lesssim \sum_{\ell
\in \Lambda} \|\mathscr{F}^{-1} \sigma_{k+\ell} \mathscr{F} \Box_k
S(t) u_0\|_{L^\infty(\mathbb{R}^n)} \nonumber\\
& \lesssim \sum_{\ell \in \Lambda} \| \sigma_{k+\ell} \mathscr{F}
\Box_k u_0\|_{L^1(\mathbb{R}^n)} \lesssim  \| \Box_k
u_0\|_{L^1(\mathbb{R}^n)}, \nonumber
\end{align}
where $\Lambda$ is as in \eqref{orthogonal}.  It follows that
\begin{align}
\|S(t)u_0\|_{M_{\infty,1}} \lesssim    \| u_0\|_{M_{1,1}}.
\label{Str-4}
\end{align}
Hence, in view of \eqref{Str-3} and \eqref{Str-4}, we have
\begin{align}
\|S(t)u_0\|_{M_{\infty,1}} \lesssim (1+ |t|)^{-n/2}  \|
u_0\|_{M_{1,1}}. \label{Str-4a}
\end{align}
By Plancherel's identity, one has that
\begin{align}
\|S(t)u_0\|_{M_{2,1}} =  \| u_0\|_{M_{2,1}}. \label{Str-5}
\end{align}
Hence, an interpolation between \eqref{Str-4a} and \eqref{Str-5}
yields (cf. \cite{WaHu}),
\begin{align}
\|S(t)u_0\|_{M_{p,1}} \lesssim (1+ |t|)^{-n(1/2-1/p)}  \|
u_0\|_{M_{p',1}}. \nonumber
\end{align}
Applying Proposition \ref{Strichartz-Mod}, we immediately obtain
that
\begin{prop}\label{Strichartz-Mod-1}
Let $2\le p<\infty$, $2/\gamma(p)= n(1/2-1/p)$.  We have for any
$\gamma\ge 2 \vee \gamma(p)$,
\begin{align}
\|S(t)u_0\|_{L^\gamma(\mathbb{R}, M^{\alpha/2}_{p,1})} \lesssim \|
u_0\|_{M_{2,1}}. \label{Str-6}
\end{align}
In particular, if $p\ge 2+4/n:=2^*$, then
\begin{align}
\|S(t)u_0\|_{L^p(\mathbb{R}, M^{\alpha/2}_{p,1})} \lesssim \|
u_0\|_{M_{2,1}}. \label{Str-7}
\end{align}
\end{prop}
Let $ \Lambda= \{\ell\in \mathbb{Z}^n:\; {\rm supp}\, \sigma_\ell
\cap {\rm supp}\,\sigma_0\not= \varnothing\}$ be as in
\eqref{orthogonal}. Using the fact that $\Box_k\Box_{k+\ell}=0$ if
$\ell \not\in \Lambda$,  it is easy to see that \eqref{Str-7}
implies the following frequency-uniform estimates:
\begin{align}
\|\,\Box_k S(t)u_0\|_{L^p_{t,x}(\mathbb{R}\times \mathbb{R}^n)}
\lesssim \|\, \Box_k u_0\|_{2}, \quad k\in \mathbb{Z}^n.
\label{Str-8}
\end{align}
Applying this estimate, we can get the following

\begin{prop} \label{prop-max-3}
Let $p\ge 2+ 4/n:= 2^*$ For any  $s>(n+2)/p$, we have
\begin{align}
\sum_{k\in \mathbb{Z}^n} \left( \sum_{\alpha\in \mathbb{Z}^n}
\|\,\Box_k S(t)u_0\|^p_{L^\infty_{t,x}(\mathbb{R} \times Q_\alpha)}
\right)^{1/p} \lesssim  \|u_0\|_{M^{s}_{2,1}}. \label{MaxMod-1}
\end{align}
\end{prop}
{\bf Proof.} Let us follow the proof of Proposition
\ref{prop-max-2}.  Denote $\langle
\partial_t \rangle=(I-\partial^2_t)^{1/2}$.  By Lemma
\ref{Max-Besov-control}, for any $s_0>1/p$,
\begin{align}
\sum_{k\in \mathbb{Z}^n} \left( \sum_{\alpha\in \mathbb{Z}^n}
\|\,\Box_k S(t)u_0\|^p_{L^\infty_{t,x}(\mathbb{R} \times Q_\alpha)}
\right)^{1/p}  & \lesssim \sum_{k\in \mathbb{Z}^n} \|\langle
\partial_t \rangle^{s_0} S(t)\Box_k u_0\|_{L^{p}(\mathbb{R},
B^{ns_0}_{p,p}(\mathbb{R}^n))} \nonumber\\
 & \lesssim \sum_{k\in \mathbb{Z}^n}
\|\langle
\partial_t \rangle^{s_0} S(t)\Box_k u_0\|_{L^{p}(\mathbb{R},
H^{ns_0}_{p}(\mathbb{R}^n))}, \label{MM-1}
\end{align}
where we have used the fact that $H^{ns_0}_{p}(\mathbb{R}^n) \subset
B^{ns_0}_{p,p}(\mathbb{R}^n)$.  Since ${\rm supp}\sigma_k \subset
B(k, \sqrt{n/2})$, applying Bernstein's multiplier estimate, we get
that
\begin{align}
\sum_{k\in \mathbb{Z}^n} \|\langle
\partial_t \rangle^{s_0} S(t)\Box_k u_0\|_{L^{p}(\mathbb{R},
H^{ns_0}_{p}(\mathbb{R}^n))} \lesssim \sum_{k\in \mathbb{Z}^n}
\langle k \rangle^{n s_0} \|\langle
\partial_t \rangle^{s_0} S(t)\Box_k u_0\|_{L^{p}_{t,x}(\mathbb{R}^{1+n})}.  \label{MM-2}
\end{align}
Similarly as in \eqref{MB-5}, using \eqref{Str-8}, we have
\begin{align}
\|\langle
\partial_t \rangle^{s_0} S(t)\Box_k
u_0\|_{L^{p}_{t,x}(\mathbb{R}^{1+n})}  & \lesssim \sum^\infty_{j=0}
\| \mathscr{F}^{-1}_{x} \langle |\xi|^2_\pm \rangle^{s_0} \varphi_j
(|\xi|^2_\pm) e^{{\rm i}t|\xi|^2_\pm} \sigma_k
(\xi) \mathscr{F}_x u_0\|_{L^{p}_{t,x}(\mathbb{R}^{1+n})} \nonumber\\
& \lesssim  \sum^\infty_{j=0} \langle k \rangle^{2 s_0}  \|
\mathscr{F}^{-1}_{x}  \varphi_j (|\xi|^2_\pm)  \sigma_k (\xi)
\mathscr{F}_x u_0\|_{L^{2}_{x}(\mathbb{R}^{n})} . \label{MM-3}
\end{align}
In an analogous way as in \eqref{MB-7} and \eqref{MB-8}, we obtain
that
\begin{align}
 \sum^\infty_{j=0} \langle k \rangle^{2 s_0}  \|
\mathscr{F}^{-1}_{x}  \varphi_j (|\xi|^2_\pm)  \sigma_k (\xi)
\mathscr{F}_x u_0\|_{L^{2}_{x}(\mathbb{R}^{n})} \lesssim \langle k
\rangle^{2 s_0+2\varepsilon}  \|\, \Box_k
u_0\|_{L^{2}_{x}(\mathbb{R}^{n})}. \label{MM-4}
\end{align}
Collecting \eqref{MM-1}--\eqref{MM-4}, we have
\begin{align}
\sum_{k\in \mathbb{Z}^n} \left( \sum_{\alpha\in \mathbb{Z}^n}
\|\,\Box_k S(t)u_0\|^p_{L^\infty_{t,x}(\mathbb{R} \times Q_\alpha)}
\right)^{1/p}  \lesssim \sum_{k\in \mathbb{Z}^n}\langle k
\rangle^{(n+2) s_0+ 2 \varepsilon} \| \,\Box_k
u_0\|_{L^{2}_{x}(\mathbb{R}^{n})}. \label{MM-5}
\end{align}
Hence, by \eqref{MM-5} we have \eqref{MaxMod-1}. $\hfill\Box$

Using the ideas as in Lemma \ref{Max-Besov-control} and Proposition
\ref{prop-max-2}, we can show the following

\begin{prop} \label{prop-max-4}
Let $p\ge 2+ 4/n:= 2^*$. Let $2^* \le r, q \le \infty$, $s_0> 1/2^*
-1/q$, $s_1 > n(1/2^* -1/r)$.  Then we have
\begin{align}
\left( \sum_{\alpha\in \mathbb{Z}^n} \|\,
S(t)u_0\|^p_{L^q(\mathbb{R},  \, L^r( Q_\alpha))} \right)^{1/p}
\lesssim \|u_0\|_{H^{s_1+2s_0}}. \label{ML-1}
\end{align}
In particular, for any $q, p\ge 2^*$,  $s> n/2- 2/q$,
\begin{align}
\left( \sum_{\alpha\in \mathbb{Z}^n} \|\,
S(t)u_0\|^p_{L^q(\mathbb{R},  \, L^\infty( Q_\alpha))} \right)^{1/p}
\lesssim \|u_0\|_{H^{s}}. \label{ML-2}
\end{align}

\end{prop}
{\bf Sketch of Proof.} In view of $\ell^{2^*} \subset \ell^p $, it
suffices to consider the case $p=2^*$. Using the inclusions
$H^{s_0}_p(\mathbb{R}) \subset L^q(\mathbb{R})$ and
$B^{s_1}_{p,p}(\mathbb{R}^n) \subset L^r(\mathbb{R}^n)$, we have
\begin{align}
 \|u\|_{L^q(\mathbb{R}, \, L^r( Q_\alpha))}
\lesssim \|(I-\partial^2_t)^{s_0/2} \sigma_\alpha
u\|_{L^p(\mathbb{R}, B^{s_1}_{p,p}(\mathbb{R}^n))}. \label{ML-4}
\end{align}
Using the same way as in
Lemma \ref{Max-Besov-control}, we can show that
\begin{align}
\left( \sum_{\alpha\in \mathbb{Z}^n} \|u\|^p_{L^q(\mathbb{R}, \,
L^r( Q_\alpha))} \right)^{1/p} \lesssim \|(I-\partial^2_t)^{s_0/2}
u\|_{L^p(\mathbb{R}, B^{s_1}_{p,p}(\mathbb{R}^n))}. \label{ML-3}
\end{align}
One can repeat the procedures as in the proof of Lemma
\ref{Max-Besov-control} to conclude that
\begin{align}
\sum_{\alpha\in \mathbb{Z}^n} \|(I-\partial^2_t)^{s_0/2}
\sigma_\alpha S(t) u_0\|^p_{L^p(\mathbb{R},
B^{s_1}_{p,p}(\mathbb{R}^n))} \lesssim \sum^\infty_{j=0}
\|\mathscr{F}^{-1}\varphi_j(|\xi|^2_\pm) \mathscr{F}
u_0\|_{H^{s_1+2s_0}(\mathbb{R}^n)}. \label{ML-5}
\end{align}
Applying an analogous way as in the proof of Proposition
\ref{prop-max-2},
\begin{align}
\sum^\infty_{j=0} \|\mathscr{F}^{-1}\varphi_j(|\xi|^2_\pm)
\mathscr{F} u_0\|_{H^{s_1+2s_0}(\mathbb{R}^n)} \lesssim
\|u_0\|_{H^{s_1+2s_0+ 2\varepsilon}}. \label{ML-6}
\end{align}
Collecting \eqref{ML-3} and \eqref{ML-6}, we immediately get
\eqref{ML-1}. $\hfill\Box$

\begin{prop} \label{prop-max-5}
For any $q \ge p\ge 2^*$,  $s> (n+2)/p- 2/q$,
\begin{align}
\sum_{k\in \mathbb{Z}^n}\left( \sum_{\alpha\in \mathbb{Z}^n}
\|\,\Box_k S(t)u_0\|^p_{L^q(\mathbb{R},  \, L^\infty( Q_\alpha))}
\right)^{1/p} \lesssim \|u_0\|_{M^{s}_{2,1}}. \label{ML-2}
\end{align}
\end{prop}

\section{Global-local estimates on time-space}
\label{global-local-estimate}

\subsection{Time-global and space-local Strichartz estimates}

We need some modifications to the Strichartz estimates, which are
global on time variable and local on spatial variable. We always
denote by $S(t)$ and $\mathscr{A}$ the generalized Schr\"odinger
semi-group and the integral operator as in \eqref{S(t)}.

\begin{prop} \label{local-Str}
Let $n\ge 3$.  Then we have
\begin{align}
\sup_{\alpha\in \mathbb{Z}^n}\|S(t) u_0\|_{L^2_{t,x}
(\mathbb{R}\times Q_\alpha)} & \lesssim \|u_0\|_2,
\label{LocStr-1}\\
\sup_{\alpha\in \mathbb{Z}^n} \left\|\mathscr{A} F
\right\|_{L^2_{t,x} (\mathbb{R}\times Q_\alpha)} & \lesssim
\sum_{\alpha\in \mathbb{Z}^n} \|F\|_{L^1_tL^2_{x} (\mathbb{R}\times
Q_\alpha)}. \label{LocStr-2}\\
\sup_{\alpha\in \mathbb{Z}^n} \left\|\mathscr{A} F
\right\|_{L^2_{t,x} (\mathbb{R}\times Q_\alpha)} & \lesssim
\sum_{\alpha\in \mathbb{Z}^n} \|F\|_{L^2_{t,x} (\mathbb{R}\times
Q_\alpha)}. \label{LocStr-2a}
\end{align}
\end{prop}
{\bf Proof.} In view of H\"older's inequality and the endpoint
Strichartz estimate,
\begin{align}
\|S(t) u_0\|_{L^2_{t,x} (\mathbb{R}\times Q_\alpha)} & \lesssim
\|S(t) u_0\|_{L^2_t L^{2n/(n-2)}_x (\mathbb{R}\times Q_\alpha)}
\nonumber\\
& \le \|S(t) u_0\|_{L^2_t L^{2n/(n-2)}_x
(\mathbb{R}\times \mathbb{R}^n)} \nonumber\\
& \lesssim \| u_0\|_{L^2_x (\mathbb{R}^n)}. \label{LocStr-3}
\end{align}
Using the above ideas and the following Strichartz estimate
\begin{align}
\left\|\mathscr{A} F \right\|_{L^2_{t}L^{2n/(n-2)}_x
(\mathbb{R}\times \mathbb{R}^n)} & \lesssim \|F\|_{L^1_tL^2_{x}
(\mathbb{R}\times \mathbb{R}^n)}, \label{LocStr-4}\\
\left\|\mathscr{A} F \right\|_{L^2_{t}L^{2n/(n-2)}_x
(\mathbb{R}\times \mathbb{R}^n)} & \lesssim
\|F\|_{L^2_tL^{2n/(n+2)}_{x} (\mathbb{R}\times \mathbb{R}^n)},
\label{LocStr-4a}
\end{align}
one can easily get \eqref{LocStr-2} and \eqref{LocStr-2a}. $\hfill
\Box$

Since the endpoint Strichartz estimates used in the proof of
Proposition \ref{local-Str} only holds for $n\ge 3$,  it is not
clear for us if \eqref{LocStr-1} still hold for $n=2$. This is why
we have an additional condition that $u_0\in \dot H^{-1/2}$ is small
in 2D. However, we have the following (see \cite{KePoVe1})

\begin{prop} \label{local-Str-2DH}
Let $n=2$. Then we have for any $1\le r<4/3$,
\begin{align}
\sup_{\alpha\in \mathbb{Z}^n} \left\|S(t) u_0 \right\|_{L^2_{t,x}
(\mathbb{R}\times Q_\alpha)} & \lesssim
\min\left(\|(-\Delta)^{-1/4}u_0\|_{2},  \; \|u_0\|_{L^2\cap
L^{r}(\mathbb{R}^n)} \right) . \label{LocStr-5}
\end{align}
\end{prop}
In the low frequency case, one easily sees that \eqref{LocStr-5} is
strictly weak than \eqref{LocStr-1}.

\medskip

\noindent{\bf Proof.} By Lemma \ref{local-sm-1}, it suffices to show
\begin{align}
\sup_{\alpha\in \mathbb{Z}^n} \left\|S(t) u_0 \right\|_{L^2_{t,x}
(\mathbb{R}\times Q_\alpha)} & \lesssim
  \|u_0\|_{L^2\cap
L^{r}(\mathbb{R}^n)}. \label{LocStr-5-1}
\end{align}
Using the unitary property in $L^2$ and the $L^p-L^{p'}$ decay
estimates of $S(t)$, we have
\begin{align}
\left\|S(t) u_0 \right\|_{L^2_{x} (Q_\alpha)}  \lesssim
(1+|t|)^{1-2/r}
  \|u_0\|_{L^2\cap
L^{r}(\mathbb{R}^n)}. \label{LocStr-5-2}
\end{align}
Taking the $L^2_t$ norm in both sides of \eqref{LocStr-5-2}, we
immediately get \eqref{LocStr-5-1}.  Hence, the result follows.
$\hfill \Box$

\begin{prop} \label{local-Str-2D}
Let $n=2$. Then we have
\begin{align}
\sup_{\alpha\in \mathbb{Z}^n} \left\|\mathscr{A} F
\right\|_{L^2_{t,x} (\mathbb{R}\times Q_\alpha)} & \lesssim
\sum_{\alpha\in \mathbb{Z}^2} \|F\|_{L^1_tL^2_{x} (\mathbb{R}\times
Q_\alpha)}. \label{LocStr-6}
\end{align}
\end{prop}
{\bf Proof.} We notice that
\begin{align}
\|S(t) f\|_{L^2_{x} (Q_\alpha)} & \lesssim  (1+|t|)^{-1}
\|f\|_{L^1_x \cap L^2_{x} (\mathbb{R}^n)}. \label{LocStr-7}
\end{align}
It follows that
\begin{align}
\ \left\|\mathscr{A} F \right\|_{L^2_{x} ( Q_\alpha)} & \lesssim
\int_{\mathbb{R}}(1+|t-\tau|)^{-1} \|F(\tau)\|_{L^1_x \cap L^2_{x}
(\mathbb{R}^n)} d\tau. \label{LocStr-8}
\end{align}
Using Young's inequality, one has that
\begin{align}
\ \left\|\mathscr{A} F \right\|_{L^2_{t, x} (\mathbb{R} \times
Q_\alpha)} & \lesssim \|F\|_{L^1(\mathbb{R},\, L^1_x \cap L^2_{x}
(\mathbb{R}^n))}. \label{LocStr-9}
\end{align}
In view of H\"older's inequality, \eqref{LocStr-9} yields the
result, as desired. $\hfill \Box$

\subsection{Note on the time-global and space-local smooth effects}

Kenig, Ponce and Vega \cite{KePoVe,KePoVe1} obtained the local
smooth effect estimates for the Schr\"odinger group $e^{{\rm i}t
\Delta}$, and their results can also be developed to the
non-elliptical Schr\"odinger group $e^{{\rm i}t \Delta_\pm}$
(\cite{KePoVe2}). On the basis of their results and Proposition
\ref{local-Str}, we can obtain a time-global version of the local
smooth effect estimates with the nonhomogeneous derivative
$(I-\Delta)^{1/2}$ instead of homogeneous derivative $\nabla$, which
is useful to control the low frequency parts of the nonlinearity.

\begin{lem} \label{local-sm-1} (\cite{KePoVe})
Let  $\Omega$ be an open set in $\mathbb{R}^n$, $\phi$ be a
$C^1(\Omega)$ function such that $\nabla \phi(\xi) \not=0$ for any
$\xi \in \Omega$. Assume that there is $N\in \mathbb{N}$ such that
for any $\bar{\xi}:=(\xi_1,...,\xi_{n-1}) \in \mathbb{R}^{n-1}$ and
$r\in \mathbb{R}$, the equation $\phi (\xi_1,...,\xi_{k}, x,
\xi_{k+1},...,\xi_{n-1})=r$ has at most $N$ solutions. For $a(x,s)
\in L^\infty(\mathbb{R}^{n}\times \mathbb{R})$ and $f\in
\mathscr{S}(\mathbb{R}^n)$,  we denote
\begin{align}
W(t)f(x) = \int_\Omega e^{{\rm i}(t\phi(\xi)+x\xi)} a(x, \phi(\xi))
\hat{f}(\xi) d\xi.  \label{LSM-1}
\end{align}
Then for $n\ge 2$, we have
\begin{align}
\|W(t)f\|_{L^2_{t,x} (\mathbb{R}\times B(0,R))}  \le CNR^{1/2}
\||\nabla \phi|^{-1/2} \hat{f}\|_{L^2(\Omega)} . \label{LSM-2}
\end{align}
\end{lem}

\begin{cor} \label{local-sm-2}
Let $n\ge 3$, $S(t)=e^{{\rm i}t \Delta_\pm}$. We have
\begin{align}
\sup_{\alpha\in \mathbb{Z}^n}\|S(t)u_0\|_{L^2_{t,x}
(\mathbb{R}\times Q_\alpha)} & \lesssim \|u_0\|_{H^{-1/2}},
\label{LSM-3}\\
 \left\|\mathscr{A} f
\right\|_{L^\infty(\mathbb{R}, H^{1/2})} & \lesssim \sum_{\alpha\in
\mathbb{Z}^n}  \|f\|_{L^2_{t,x} (\mathbb{R}\times Q_\alpha)}.
\label{LSM-3-dual}
\end{align}
For $n=2$, \eqref{LSM-3-dual} also holds if one substitutes
$H^{1/2}$ by $\dot H^{1/2}$.
\end{cor}
{\bf Proof.} Let $\Omega=\mathbb{R}^n\setminus B(0,1)$, $\phi(\xi)=
|\xi|^2_\pm$ and $\psi$ be as in \eqref{cutoff}, $a(x, s)=
1-\psi(s)$ in Lemma \ref{local-sm-1}. Taking $W(t):=
S(t)\mathscr{F}^{-1}(1-\psi)\mathscr{F}$, from \eqref{LSM-2} we have
\begin{align}
\sup_{\alpha \in \mathbb{Z}^n}
\|S(t)\mathscr{F}^{-1}(1-\psi)\mathscr{F} u_0\|_{L^2_{t,x}
(\mathbb{R}\times Q_\alpha)} \lesssim
\||\xi|^{-1/2}\hat{u}_0\|_{L^2_\xi(\mathbb{R}^n\setminus B(0,1))}.
\label{LSM-4}
\end{align}
It follows from Proposition \ref{local-Str} that
\begin{align}
\|S(t)\mathscr{F}^{-1}\psi\mathscr{F} u_0\|_{L^2_{t,x}
(\mathbb{R}\times Q_\alpha)} & \lesssim
\|\mathscr{F}^{-1}\psi\mathscr{F} u_0\|_{L^2_x
(\mathbb{R}^n)}\nonumber\\
&  \lesssim \|\hat{u}_0\|_{L^2_\xi (B(0,2))}. \label{LSM-6}
\end{align}
From \eqref{LSM-4} and \eqref{LSM-6} we have \eqref{LSM-3}, as
desired.  \eqref{LSM-3-dual} is the dual version of \eqref{LSM-3}.
$\hfill \Box$

\medskip

When $n=2$,  it is known that for the elliptic case, the endpoint
Strichartz estimate holds for the radial function (cf. \cite{Tao}).
So, Corollary \ref{local-sm-2} also holds for the radial function
$u_0$ in the elliptic case. The following local smooth effect
estimates for the nonhomogeneous part of the solutions of the
Schr\"odinger equation is also due to Kenig, Ponce and Vega
\cite{KePoVe1}\footnote{In \cite{KePoVe1}, the result was stated for
the elliptic case, however, their result is also adapted to the
non-elliptic cases.}.

\begin{prop} \label{local-sm-3}
Let $n\ge 2$, $S(t)=e^{{\rm i}t \Delta_\pm}$. We have
\begin{align}
\sup_{\alpha\in \mathbb{Z}^n} \left\| \nabla \mathscr{A} f
\right\|_{L^2_{t,x} (\mathbb{R}\times Q_\alpha)} & \lesssim
\sum_{\alpha\in \mathbb{Z}^n} \|f\|_{L^2_{t,x} (\mathbb{R}\times
Q_\alpha)}. \label{LSM-7}
\end{align}
\end{prop}

\section{Proof of Theorem \ref{GWP-nD}} \label{proof-gwp-nd}

\begin{lem} \label{Sobolev-ineq}
(Sobolev Inequality).  Let $\Omega \subset \mathbb{R}^n$ be a
bounded domain with $\partial \Omega \in C^m$, $m, \ell \in
\mathbb{N}\cup \{0\}$, $1\le r,p,q\le \infty$. Assume that
$$
\frac{\ell}{m} \le \theta \le 1, \ \ \frac{1}{p} -\frac{\ell}{n} =
\theta\left(\frac{1}{r} - \frac{m}{n}\right) +\frac{1-\theta}{q}.
$$
Then we have
\begin{align}
\sum_{|\beta|=\ell}\|D^\beta u\|_{L^p(\Omega)} \lesssim  \|
u\|^{1-\theta}_{L^q(\Omega)} \|u\|^\theta_{W^m_r(\Omega)},
\label{Sob-1}
\end{align}
where $\|u\|_{W^m_r(\Omega)} = \sum_{|\beta|\le m}  \|
u\|_{L^r(\Omega)}$.
\end{lem}
{\bf Proof of Theorem \ref{gNLS}.} In order to illustrate our ideas
in an exact way, we first consider a simple case $s= [n/2]+ 5/2$ and
there is no difficulty to generalize the proof to the case
$s>n/2+3/2$, $s+1/2 \in \mathbb{N}$. We assume without loss of
generality that
\begin{align}
 F(u,\, \bar{u},\,
\nabla u,\, \nabla \bar{u}):= F(u, \nabla u)= \sum_{\Lambda_{\kappa,
\nu}} c_{\kappa\nu_1...\nu_n} u^{\kappa}
u^{\nu_1}_{x_1}...u^{\nu_n}_{x_n}, \label{poly-1}
\end{align}
where
$$
\Lambda_{\kappa, \nu}=\{(\kappa, \nu_1,...,\nu_n): \, m+1\le
\kappa+\nu_1+...+\nu_n \le M+1 \}.
$$
Since we only use the Sobolev norm to control the nonlinear terms,
$\bar{u}$ and $u$ have the same norm, whence, the general cases can
be handled in the same way. Denote
\begin{align}
& \lambda_1 (v) : = \|v\|_{\ell^\infty_\alpha (L^2_{t,x}(\mathbb{R}
\times Q_\alpha))}, \nonumber\\
& \lambda_2 (v) : =  \|v\|_{\ell^{2^*}_\alpha (L^\infty_{t,x}
(\mathbb{R} \times Q_\alpha))}, \nonumber\\
 & \lambda_3 (v) : =
\|v\|_{\ell^{2^*}_\alpha (L^{2m}_{t}L^\infty_ x(\mathbb{R} \times
Q_\alpha))}. \nonumber
\end{align}
Put
\begin{align}
& \mathscr{D}_n= \left\{ u: \; \sum_{|\beta|\le [n/2]+3}  \lambda_1
(D^\beta u)  + \sum_{|\beta|\le 1} \sum_{i=2,3} \lambda_i (D^\beta
u) \le \varrho  \right\}. \label{metricspace}
\end{align}
We consider the mapping
\begin{align}
& \mathscr{T}:  u(t) \to  S(t)u_0 - {\rm i} \mathscr{A} F(u, \,
\nabla u), \label{map}
\end{align}
and we show that $\mathscr{T}: \mathscr{D}_n\to \mathscr{D}_n$ is a
contraction mapping for any $n\ge 2$.

{\it Step} 1. For any $u\in \mathscr{D}_n$, we estimate $\lambda_1
(D^\beta \mathscr{T}u)$, $|\beta| \le 3+[n/2]$. We consider the
following three cases.

{\it Case } 1. $n\ge 3$ and $1\le |\beta| \le 3+[n/2]$. In view of
Corollary \ref{local-sm-2} and Proposition \ref{local-sm-3}, we have
for any $\beta$, $1\le |\beta|\le 3+[n/2]$,
\begin{align}
\lambda_1 (D^\beta \mathscr{T} u)  & \lesssim \|S(t) D^\beta
u_0\|_{\ell^\infty_\alpha (L^2_{t,x}(\mathbb{R} \times Q_\alpha))} +
\sum_{\Lambda_{\kappa, \nu}}
 \|\mathscr{A} D^\beta(u^{\kappa}
u^{\nu_1}_{x_1}...u^{\nu_n}_{x_n})\|_{\ell^\infty_\alpha
(L^2_{t,x}(\mathbb{R} \times Q_\alpha))} \nonumber\\
& \lesssim \|u_0\|_{H^s }  +  \sum_{ |\beta|\le 2+[n/2]}
\sum_{\Lambda_{\kappa, \nu}} \sum_{\alpha\in \mathbb{Z}^n}
\|D^\beta(u^{\kappa}
u^{\nu_1}_{x_1}...u^{\nu_n}_{x_n})\|_{L^2_{t,x}(\mathbb{R} \times
Q_\alpha)}. \label{Lambda1-1}
\end{align}
For simplicity, we can further assume that $u^{\kappa}
u^{\nu_1}_{x_1}...u^{\nu_n}_{x_n} = u^{\kappa} u^{\nu}_{x_1}$ in
\eqref{Lambda1-1} and the general case can be treated in an
analogous way\footnote{One can see below for a general treating.}.
So, one can rewrite \eqref{Lambda1-1} as
\begin{align}
\sum_{ 1\le |\beta|\le 3+[n/2]}\lambda_1 (D^\beta \mathscr{T} u)
\lesssim \|u_0\|_{H^s }  +  \sum_{ |\beta|\le 2+[n/2]}
\sum_{\Lambda_{\kappa, \nu}} \sum_{\alpha\in \mathbb{Z}^n}
\|D^\beta(u^{\kappa} u^{\nu}_{x_1})\|_{L^2_{t,x}(\mathbb{R} \times
Q_\alpha)}. \label{Lambda1-2}
\end{align}
It is easy to see that
\begin{align}
|D^\beta(u^{\kappa} u^{\nu}_{x_1})|  \lesssim
\sum_{\beta_1+...+\beta_{\kappa+\nu}=\beta} | D^{\beta_1}
u....D^{\beta_\kappa} u D^{\beta_{\kappa+1}}
u_{x_1}...D^{\beta_{\kappa+\nu}} u_{x_1} |. \label{Lambda1-3}
\end{align}
By H\"older's inequality,
\begin{align}
\|D^\beta(u^{\kappa} u^{\nu}_{x_1})\|_{L^2_x(Q_\alpha)} \lesssim
\sum_{\beta_1+...+\beta_{\kappa+\nu}=\beta} \prod^\kappa_{i=1} \|
D^{\beta_i}u\|_{L^{p_i}_x(Q_\alpha)} \prod^{\kappa+\nu}_{i=\kappa+1}
\| D^{\beta_i} u_{x_1}\|_{L^{p_i}_x(Q_\alpha)}, \label{Lambda1-4}
\end{align}
where
$$
p_i=\left\{
\begin{array}{ll}
2|\beta|/|\beta_i|, & |\beta_i|\ge 1,\\
\infty, & |\beta_i|=0.
\end{array}
\right.
$$
It is easy to see that for $\theta_i= |\beta_i|/|\beta|$,
$$
\frac{1}{p_i}- \frac{|\beta_i|}{n} = \theta_i \left( \frac{1}{2}-
\frac{|\beta|}{n}\right)+ \frac{1-\theta_i}{\infty}.
$$
Using Sobolev's inequality, one has that for $B_\alpha :=\{x:\;
|x-\alpha|\le \sqrt{n}\}$,
\begin{align}
& \|D^{\beta_i}u\|_{L^{p_i}_x(Q_\alpha)} \le
\|D^{\beta_i}u\|_{L^{p_i}_x(B_\alpha)} \lesssim
\|u\|^{1-\theta_i}_{L^{\infty}_x(B_\alpha)}
 \|u \|^{\theta_i}_{W^{|\beta|}_2(B_\alpha)}, \ \ i=1,...,\kappa;
 \label{Lambda1-5}\\
& \|D^{\beta_i}u_{x_1}\|_{L^{p_i}_x(Q_\alpha)} \lesssim
\|u_{x_1}\|^{1-\theta_i}_{L^{\infty}_x(B_\alpha)}
 \|u_{x_1} \|^{\theta_i}_{W^{|\beta|}_2(B_\alpha)}, \ \ i=\kappa+1,..., \kappa+\nu. \label{Lambda1-6}
\end{align}
Since
$$
\sum^{\kappa+\nu}_{i=1} \theta_i =1, \quad  \sum^{\kappa+\nu}_{i=1}
(1-\theta_i) = \kappa+\nu-1,
$$
by \eqref{Lambda1-4}--\eqref{Lambda1-6} we have
\begin{align}
\|D^\beta(u^{\kappa} u^{\nu}_{x_1})\|_{L^2_x(Q_\alpha)} \lesssim &
\sum_{|\beta| \le 2+[n/2]} (\|u \|_{W^{|\beta|}_2(B_\alpha)}+
\|u_{x_1}\|_{W^{|\beta|}_2(B_\alpha)}) \nonumber\\
& \times (\|u\|^{\kappa+\nu-1}_{L^{\infty}_x(B_\alpha)}+
\|u_{x_1}\|^{\kappa+\nu-1}_{L^{\infty}_x(B_\alpha)}) \nonumber\\
\lesssim & \sum_{|\gamma| \le 3+[n/2]} \|D^\gamma u \|_{L^2_x
(B_\alpha)} \sum_{|\beta|\le 1}\|D^\beta
u\|^{\kappa+\nu-1}_{L^{\infty}_x(B_\alpha)}. \label{Lambda1-7}
\end{align}
It follows from \eqref{Lambda1-7} and $\ell^{2^*} \subset
\ell^{\kappa+\nu-1}$ that
\begin{align}
 & \!\!\!\! \sum  _{|\beta| \le 2+[n/2]} \sum_{\alpha\in \mathbb{Z}^n}
\|D^\beta(u^{\kappa} u^{\nu}_{x_1})\|_{L^2_{t,x}(\mathbb{R}\times Q_\alpha)} \nonumber\\
& \lesssim   \sum_{\alpha\in \mathbb{Z}^n} \sum_{|\gamma| \le
3+[n/2]} \|D^\gamma u \|_{L^2_{t,x} (\mathbb{R}\times B_\alpha)}
\sum_{|\beta|\le 1}\|D^\beta
u\|^{\kappa+\nu-1}_{L^{\infty}_{t,x}(\mathbb{R}\times B_\alpha)}
\nonumber\\
& \lesssim    \sum_{|\gamma| \le 3+[n/2]} \lambda_1(D^\gamma u)
 \sum_{|\beta|\le
1}\lambda_2 (D^\beta u)^{\kappa+\nu-1} \lesssim
\varrho^{\kappa+\nu}. \label{Lambda1-8}
\end{align}
Hence, in view of \eqref{Lambda1-2} and \eqref{Lambda1-8} we have
\begin{align}
\sum_{ 1\le |\beta|\le 3+[n/2]}\lambda_1 (D^\beta \mathscr{T} u)
\lesssim \|u_0\|_{H^s }  + \sum^{M+1}_{\kappa+\nu=m+1}
\varrho^{\kappa+\nu}. \label{Lambda1-9}
\end{align}

{\it Case } 2. $n\ge 3$ and $|\beta|=0$. By Corollary
\ref{local-sm-2}, the local Strichartz estimate \eqref{LocStr-2}
 and H\"older's inequality,
\begin{align}
\lambda_1 (\mathscr{T}u)  & \lesssim \|S(t)
u_0\|_{\ell^\infty_\alpha (L^2_{t,x}(\mathbb{R} \times Q_\alpha))} +
 \|\mathscr{A} F(u, \nabla u)\|_{\ell^\infty_\alpha
(L^2_{t,x}(\mathbb{R} \times Q_\alpha))} \nonumber\\
& \lesssim \|u_0\|_{2}  +  \sum_{\alpha\in \mathbb{Z}^n} \|F(u,
\nabla u) \|_{L^1_{t}L^{2}_x(\mathbb{R} \times
Q_\alpha)} \nonumber\\
& \lesssim \|u_0\|_{2}  +  \sum_{\Lambda_{\kappa, \nu}}
\sum_{\alpha\in \mathbb{Z}^n} \|u^\kappa u^{\nu_1}_{x_1}...
u^{\nu_n}_{x_n} \|_{L^1_{t}L^{2}_x(\mathbb{R} \times
Q_\alpha)} \nonumber\\
& \lesssim \|u_0\|_{2}  +  \sum^{M+1}_{\kappa+\nu=m+1}
\sum_{|\gamma|\le 1} \sup_{\alpha\in \mathbb{Z}^n}   \|D^\gamma u
\|_{L^2_{t,x} (\mathbb{R} \times Q_\alpha)} \nonumber\\
& \ \ \ \ \ \ \ \ \ \ \ \times  \sum_{|\beta|\le 1} \sum_{\alpha\in
\mathbb{Z}^n} \|D^\beta u \|^{\kappa+\nu-1}_{L^{2(\kappa+\nu-1)}_{t}
L^\infty_x (\mathbb{R} \times
Q_\alpha)}  \nonumber\\
& \lesssim \|u_0\|_{2}  +  \sum^{M+1}_{\kappa+\nu=m+1}
\sum_{|\gamma|\le 1} \lambda_1(D^\gamma u) \sum_{i=2,3}
\sum_{|\beta|\le 1} \lambda_i(D^\beta u )^{\kappa+\nu-1}
\nonumber\\
& \lesssim \|u_0\|_{2}  +  \sum^{M+1}_{\kappa+\nu=m+1}
 \varrho^{\kappa+\nu}. \label{Lambda1-10}
\end{align}
{\it Case } 3. $n=2, \; |\beta|=0$. By Propositions
\ref{local-Str-2DH} and \ref{local-Str-2D}, we have
\begin{align}
\lambda_1 (\mathscr{T}u) & \lesssim \|u_0\|_{\dot H^{-1/2}} +
\sum_{\alpha\in \mathbb{Z}^n} \|F(u, \nabla u)
\|_{L^1_{t}L^{2}_x(\mathbb{R} \times Q_\alpha)}. \label{Lambda1-11}
\end{align}
Using the same way as in Case 2, we have
\begin{align}
\lambda_1 (\mathscr{T}u) & \lesssim \|u_0\|_{\dot H^{-1/2}} +
\sum^{M+1}_{\kappa+\nu=m+1} \varrho^{\kappa+\nu}. \label{Lambda1-12}
\end{align}

 {\it Step} 2. We consider the estimates of $\lambda_2(D^\beta
 \mathscr{T}u)$, $|\beta| \le 1$. Using the estimates of the maximal
 function as in Proposition \ref{prop-max-2},  we have for $|\beta| \le
 1$, $0<\varepsilon \ll 1$,
\begin{align}
\lambda_2 (D^\beta \mathscr{T}u)  & \lesssim \|S(t)D^\beta
u_0\|_{\ell^{2^*}_\alpha (L^\infty_{t,x}(\mathbb{R} \times
Q_\alpha))} +
 \|\mathscr{A} D^\beta F(u, \nabla u)\|_{\ell^{2^*}_\alpha
(L^\infty_{t,x}(\mathbb{R} \times Q_\alpha))} \nonumber\\
& \lesssim \|D^\beta u_0\|_{H^{n/2+\varepsilon}}   +  \sum_{|\beta|
\le 1} \|D^\beta F(u, \nabla u)
\|_{L^1 (\mathbb{R}, H^{n/2+\varepsilon}(\mathbb{R}^n))} \nonumber\\
& \lesssim   \|u_0\|_{H^{n/2+1+\varepsilon}} + \sum_{|\beta| \le
[n/2]+2}\sum_{\alpha\in \mathbb{Z}^n} \|D^\beta F(u, \nabla u)
\|_{L^1_{t}L^{2}_x(\mathbb{R} \times Q_\alpha)}. \label{Lambda2-1}
\end{align}
Applying the same way as in Step 1, for any $|\beta| \le [n/2]+2$,
\begin{align}
\|D^\beta F(u, \nabla u)\|_{L^2_x(Q_\alpha)}  \lesssim
\sum^{M+1}_{\kappa+\nu=m+1} \sum_{|\beta|\le 1}\|D^\beta
u\|^{\kappa+\nu-1}_{L^{\infty}_x(B_\alpha)}  \sum_{|\gamma| \le
3+[n/2]} \|D^\gamma u \|_{L^2_x (B_\alpha)} . \label{Lambda2-2}
\end{align}
By H\"older's inequality, we have from \eqref{Lambda2-2} that
\begin{align}
\|D^\beta F(u, \nabla u)\|_{L^1_tL^2_x(\mathbb{R}\times Q_\alpha)}
\lesssim & \sum^{M+1}_{\kappa+\nu=m+1}   \sum_{|\gamma| \le 3+[n/2]}
\|D^\gamma u \|_{L^2_{t,x} (\mathbb{R}\times
B_\alpha)} \nonumber\\
&  \ \ \ \ \ \times \sum_{|\beta|\le 1}\|D^\beta
u\|^{\kappa+\nu-1}_{L^{2(\kappa+\nu-1)}_t
L^{\infty}_x(\mathbb{R}\times B_\alpha)}. \label{Lambda2-3}
\end{align}
Summarizing \eqref{Lambda2-3} over all $\alpha\in \mathbb{Z}^n$, we
have for any $|\beta| \le 2+[n/2]$,
\begin{align}
& \sum_{\alpha \in \mathbb{Z}^n} \|D^\beta F(u, \nabla
u)\|_{L^1_tL^2_x(\mathbb{R}\times Q_\alpha)} \nonumber\\
 & \lesssim
\sum^{M+1}_{\kappa+\nu=m+1}    \sum_{|\gamma| \le 3+[n/2]} \lambda_1
( D^\gamma u )  \sum_{|\beta|\le 1} \sum_{\alpha \in \mathbb{Z}^n}
\|D^\beta u\|^{\kappa+\nu-1}_{L^{2(\kappa+\nu-1)}_t
L^{\infty}_x(\mathbb{R}\times B_\alpha)} \nonumber\\
& \lesssim \sum^{M+1}_{\kappa+\nu=m+1}    \sum_{|\gamma| \le
3+[n/2]} \lambda_1 ( D^\gamma u )  \sum_{|\beta|\le 1} \sum_{\alpha
\in \mathbb{Z}^n} \|D^\beta u\|^{\kappa+\nu-1}_{(L^{2m}_t
L^{\infty}_x) \cap
L^\infty_{t,x}(\mathbb{R}\times B_\alpha)} \nonumber\\
& \lesssim \sum^{M+1}_{\kappa+\nu=m+1}    \sum_{|\gamma| \le
3+[n/2]} \lambda_1 ( D^\gamma u )  \sum_{|\beta|\le 1}
(\lambda_2(D^\beta u)^{\kappa+\nu-1}
+\lambda_3(D^\beta u)^{\kappa+\nu-1} ) \nonumber\\
 & \lesssim
\sum^{M+1}_{\kappa+\nu=m+1}    \varrho^{\kappa+\nu}.
\label{Lambda2-4}
\end{align}
Combining \eqref{Lambda2-1} with  \eqref{Lambda2-4}, we obtain that
\begin{align}
\sum_{|\beta|\le 1}\lambda_2 (D^\beta \mathscr{T}u)  \lesssim
\|u_0\|_{H^{n/2+1+\varepsilon}} + \sum^{M+1}_{\kappa+\nu=m+1}
\varrho^{\kappa+\nu}.  \label{Lambda2-5}
\end{align}

{\it Step} 3. We estimate $\lambda_3 (D^\beta
 \mathscr{T}u)$, $|\beta| \le 1$. In view of Proposition
 \ref{prop-max-4}, one has that
\begin{align}
\lambda_3 (D^\beta \mathscr{T}u)  & \lesssim \|S(t)D^\beta
u_0\|_{\ell^{2^*}_\alpha (L^{2m}_t L^\infty_{x}(\mathbb{R} \times
Q_\alpha))} +
 \|\mathscr{A} D^\beta F(u, \nabla u)\|_{\ell^{2^*}_\alpha
(L^{2m}_t L^\infty_{x}(\mathbb{R} \times Q_\alpha))} \nonumber\\
& \lesssim \|D^\beta u_0\|_{H^{n/2 -1/m +\varepsilon}}   +
\sum_{|\beta| \le 1} \|D^\beta F(u, \nabla u)
\|_{L^1 (\mathbb{R}, H^{n/2-1/m +\varepsilon}(\mathbb{R}^n))} \nonumber\\
& \lesssim   \|u_0\|_{H^{n/2+1}} + \sum_{|\beta| \le
[n/2+2]}\sum_{\alpha\in \mathbb{Z}^n} \|D^\beta F(u, \nabla u)
\|_{L^1_{t}L^{2}_x(\mathbb{R} \times Q_\alpha)}, \label{Lambda3-1}
\end{align}
which reduces to the case as in \eqref{Lambda2-1}.

Therefore, collecting the estimates as in Steps 1--3, we have for
$n\ge 3$,
\begin{align}
\sum_{|\beta|\le 3+[n/2]}\lambda_1 (D^\beta \mathscr{T}u) +
\sum_{i=2,3} \sum_{|\beta|\le 1}\lambda_i (D^\beta \mathscr{T}u)
\lesssim \|u_0\|_{H^s} + \sum^{M+1}_{\kappa+\nu=m+1}
\varrho^{\kappa+\nu},  \label{Lambda-1}
\end{align}
and for $n\ge 2$,
\begin{align}
\sum_{|\beta|\le 4}\lambda_1 (D^\beta \mathscr{T}u) + \sum_{i=2,3}
\sum_{|\beta|\le 1}\lambda_i (D^\beta \mathscr{T}u) \lesssim
\|u_0\|_{H^{7/2} \cap \dot H^{-1/2}} + \sum^{M+1}_{\kappa+\nu=m+1}
\varrho^{\kappa+\nu}. \label{Lambda-2}
\end{align}
It follows that for $n\ge 3,$ $\mathscr{T}: \mathscr{D}_n \to
\mathscr{D}_n$ is a contraction mapping if $\varrho$ and
$\|u_0\|_{H^s}$
 are small enough (similarly for $n=2$).

Before considering the case $s>n/2+3/2$, we first establish a
nonlinear mapping estimate:
\begin{lem} \label{NLE-Besov}
Let $n\ge 2$, $s>0$, $K\in \mathbb{N}$. Let $1\le p, p_i, q, q_i \le
\infty$ satisfy $1/p = 1/p_1+(K-1)/p_2$ and $1/q=1/q_1+ (K-1)/q_2$.
We have
\begin{align}
 \|v_1...v_K\|_{\ell^{1,s}_\triangle \ell^1_\alpha (L^q_t L^p_x
(\mathbb{R}\times Q_\alpha))}  \lesssim &
 \sum^K_{k=1}  \|v_k\|_{\ell^{1,s}_\triangle \ell^\infty_\alpha (L^{q_1}_t L^{p_1}_x
(\mathbb{R}\times Q_\alpha))} \nonumber\\
&  \  \times \prod_{i\not= k, \, i=1,...,K}
\|v_i\|_{\ell^{1}_\triangle \ell^{K-1}_\alpha (L^{q_2}_t L^{p_2}_x
(\mathbb{R}\times Q_\alpha))}. \label{NLE-1}
\end{align}
\end{lem}
{\bf Proof.} Denote $S_r u= \sum_{j\le r} \triangle u$. We have
\begin{align}
v_1... v_K = \sum^\infty_{r=-1} (S_{r+1}v_1 ... S_{r+1} v_K -
S_{r}v_1 ... S_r v_K), \label{NLE-2}
\end{align}
where we assume that $S_{-1}v \equiv 0$. Recall the identity,
\begin{align}
\prod^K_{k=1} a_k - \prod^K_{k=1} b_k = \sum^K_{k=1} (a_k-b_k)
\prod_{i\le k-1} b_i \prod_{i\ge k+1} a_i, \label{NLE-3}
\end{align}
where we assume that $\prod_{i\le 0}a_i= \prod_{i\ge K+1} \equiv 1$.
We have
\begin{align}
v_1... v_K = \sum^\infty_{r=-1} \sum^K_{k=1} \triangle_{r+1} v_k
\prod^{k-1}_{i=1} S_{r}v_i  \prod^{K}_{i=k+1} S_{r+1}v_i.
\label{NLE-4}
\end{align}
Hence, it follows that
\begin{align}
&  \|v_1...v_K\|_{\ell^{1,s}_\triangle \ell^1_\alpha (L^q_t L^p_x
(\mathbb{R}\times Q_\alpha))} \nonumber\\
& =  \sum^\infty_{j=0} 2^{sj} \sum_{\alpha\in \mathbb{Z}^n}
\|\triangle_j (v_1...v_K)\|_{L^q_t L^p_x (\mathbb{R}\times
Q_\alpha)} \nonumber\\
 &
\lesssim  \sum^K_{k=1} \sum^\infty_{j=0} 2^{sj} \sum_{\alpha\in
\mathbb{Z}^n}  \sum^\infty_{r=-1}  \left\|\triangle_j (
\triangle_{r+1} v_k \prod^{k-1}_{i=1} S_{r}v_i  \prod^{K}_{i=k+1}
S_{r+1}v_i )\right\|_{L^q_t L^p_x (\mathbb{R}\times Q_\alpha)}.
\label{NLE-5}
\end{align}
Using the support property of $\widehat{\triangle_r v}$ and
$\widehat{S_r v}$, we see that
\begin{align}
\triangle_j ( \triangle_{r+1} v_k \prod^{k-1}_{i=1} S_{r}v_i
\prod^{K}_{i=k+1} S_{r+1}v_i ) \equiv 0, \ \ j> r+C. \label{NLE-6}
\end{align}
Using the fact $\|\int f d \mu\|_X \le \int \|f\|_X d\mu$, one has
that
\begin{align}
\sum_{\alpha\in \mathbb{Z}^n}\|\triangle_j f\|_{L^q_t L^p_x
(\mathbb{R}\times Q_\alpha)} & \le \sum_{\alpha\in \mathbb{Z}^n}
\int_{\mathbb{R}^n} |\mathscr{F}^{-1}\varphi_j (y)| \|f(t,
x-y)\|_{L^q_t L^p_x (\mathbb{R}\times Q_\alpha)} dy \nonumber\\
& \le \sup_{y \in \mathbb{R}^n} \sum_{\alpha\in \mathbb{Z}^n} \|f(t,
x-y)\|_{L^q_t L^p_x (\mathbb{R}\times Q_\alpha)} \int_{\mathbb{R}^n}
|\mathscr{F}^{-1}\varphi_j (y)|  dy \nonumber\\
& \lesssim  \sum_{\alpha\in \mathbb{Z}^n} \|f\|_{L^q_t L^p_x
(\mathbb{R}\times Q_\alpha)}. \label{NLE-7}
\end{align}
Collecting \eqref{NLE-5}--\eqref{NLE-7} and using Fubini's Theorem,
we have
\begin{align}
&  \|v_1...v_K\|_{\ell^{1,s}_\triangle \ell^1_\alpha (L^q_t L^p_x
(\mathbb{R}\times Q_\alpha))} \nonumber\\
 &
\lesssim  \sum^K_{k=1}  \sum^\infty_{r=-1} \sum_{j\le r+C} 2^{sj}
\sum_{\alpha\in \mathbb{Z}^n}   \left\|\triangle_j ( \triangle_{r+1}
v_k \prod^{k-1}_{i=1} S_{r}v_i  \prod^{K}_{i=k+1} S_{r+1}v_i
)\right\|_{L^q_t L^p_x (\mathbb{R}\times Q_\alpha)} \nonumber\\
 &
\lesssim  \sum^K_{k=1}  \sum^\infty_{r=-1} \sum_{j\le r+C} 2^{sj}
 \sum_{\alpha\in \mathbb{Z}^n}
\left\|\triangle_{r+1} v_k \prod^{k-1}_{i=1} S_{r}v_i
\prod^{K}_{i=k+1} S_{r+1}v_i
\right\|_{L^q_t L^p_x (\mathbb{R}\times Q_\alpha)} \nonumber\\
 &
\lesssim  \sum^K_{k=1}  \sum^\infty_{r=-1}  2^{s r}  \sum_{\alpha\in
\mathbb{Z}^n} \left\| \triangle_{r+1} v_k \prod^{k-1}_{i=1} S_{r}v_i
\prod^{K}_{i=k+1} S_{r+1}v_i
\right\|_{L^q_t L^p_x (\mathbb{R}\times Q_\alpha)} \nonumber\\
 &
\lesssim  \sum^K_{k=1}  \sum^\infty_{r=-1}  2^{s r}  \sum_{\alpha\in
\mathbb{Z}^n} \|\triangle_{r+1} v_k\|_{L^{q_1}_t L^{p_1}_x
(\mathbb{R}\times Q_\alpha)}  \prod_{i\not= k, \ i=1,...,K}
\|v_i\|_{\ell^1_\triangle (L^{q_2}_t
   L^{p_2}_x (\mathbb{R}\times
Q_\alpha))}    \nonumber\\
 &
\lesssim  \sum^K_{k=1} \|v_k\|_{\ell^{1,s}_\triangle
\ell^\infty_\alpha  (L^{q_1}_t L^{p_1}_x (\mathbb{R}\times
Q_\alpha))} \sum_{\alpha\in \mathbb{Z}^n}
   \prod_{i\not= k, \ i=1,...,K} \|v_i\|_{\ell^1_\triangle (L^{q_2}_t
   L^{p_2}_x (\mathbb{R}\times
Q_\alpha) )}, \label{NLE-8}
\end{align}
the result follows. $\hfill\Box$

For short, we write $\|\nabla u\|_X = \|\partial_{x_1} u\|_X+...+
\|\partial_{x_n} u\|_X$.

\begin{lem} \label{LE-Besov-1}
Let $n \ge 3$. We have for any $s>0$
\begin{align}
\sum_{k=0,1}\| S(t) \nabla^k u_0\|_{\ell^{1,s}_\triangle
\ell^\infty_\alpha (L^2_{t, x} (\mathbb{R}\times Q_\alpha))} &
\lesssim \|u_0\|_{B^{s+1/2}_{2,1}},  \label{Besov-LE-1}\\
 \sum_{k=0,1}\|
\mathscr{A} \nabla^k F\|_{\ell^{1,s}_\triangle \ell^\infty_\alpha
(L^2_{t, x} (\mathbb{R}\times Q_\alpha))} & \lesssim
\|F\|_{\ell^{1,s}_\triangle \ell^1_\alpha (L^2_{t, x}
(\mathbb{R}\times Q_\alpha))}.  \label{Besov-LE-2}
\end{align}
\end{lem}
{\bf Proof.} In view of Corollary \ref{local-sm-2} and Propositions
\ref{local-Str} and \ref{local-sm-3}, we have the results, as
desired. $\hfill \Box$

\begin{lem} \label{LE-Besov-2}
Let $n =2$. We have for any $s>0$
\begin{align}
\sum_{k=0,1}\| S(t) \nabla^k u_0\|_{\ell^{1,s}_\triangle
\ell^\infty_\alpha (L^2_{t, x} (\mathbb{R}\times Q_\alpha))} &
\lesssim \|u_0\|_{B^{s+1/2}_{2,1} \cap \dot H^{-1/2}},  \label{Besov-LE-3}\\
\| \mathscr{A} \nabla F\|_{\ell^{1,s}_\triangle \ell^\infty_\alpha
(L^2_{t, x} (\mathbb{R}\times Q_\alpha))} & \lesssim
\|F\|_{\ell^{1,s}_\triangle \ell^1_\alpha (L^2_{t, x}
(\mathbb{R}\times Q_\alpha))},   \label{Besov-LE-4}\\
\|\mathscr{A}  F\|_{\ell^{1,s}_\triangle \ell^\infty_\alpha (L^2_{t,
x} (\mathbb{R}\times Q_\alpha))} & \lesssim
\|F\|_{\ell^{1,s}_\triangle \ell^1_\alpha (L^1_t L^2_{x}
(\mathbb{R}\times Q_\alpha))}.  \label{Besov-LE-5}
\end{align}
\end{lem}
{\bf Proof.} By Propositions \ref{local-Str-2DH}, \ref{local-Str-2D}
and \ref{local-sm-3}, we have the results, as desired. $\hfill \Box$

We now continue the proof of Theorem \ref{GWP-nD} and now we
consider the general case $s>n/2+3/2$. We write
\begin{align}
& \lambda_1 (v) : = \sum_{i=0,1} \|\nabla^i
v\|_{\ell^{1,s-1/2}_\triangle \ell^\infty_\alpha
(L^2_{t,x}(\mathbb{R}
\times Q_\alpha))}, \nonumber\\
& \lambda_2 (v) : = \sum_{i=0,1} \|\nabla^i
v\|_{\ell^{1,s-1/2}_\triangle \ell^{2^*}_\alpha (L^\infty_{t,x}
(\mathbb{R} \times Q_\alpha))}, \nonumber\\
 & \lambda_3 (v) : = \sum_{i=0,1}
\|\nabla^i v\|_{\ell^{1,s-1/2}_\triangle \ell^{2^*}_\alpha
(L^{2m}_{t}L^\infty_ x(\mathbb{R} \times Q_\alpha))}, \nonumber\\
& \mathscr{D} = \{u: \sum_{i=1,2,3}\lambda_i (v) \le \varrho \}.
\label{spaceD}
\end{align}
Note $\lambda_i$ and $\mathscr{D}$ defined here are different from
those in the above. We only give the details of the proof in the
case $n\ge 3$ and the case $n=2$ can be shown in a slightly
different way. Let $\mathscr{T}$ be defined as in \eqref{map}. Using
Lemma \ref{LE-Besov-1}, we have
\begin{align}
& \lambda_1 ( \mathscr{T} u)  \lesssim  \|u_0\|_{B^{s}_{2,1}} +\|F\|
_{\ell^{1,s-1/2}_\triangle \ell^1_\alpha (L^2_{t,x}(\mathbb{R}
\times Q_\alpha))}.  \label{proof-thm-1}
\end{align}
For simplicity, we write
\begin{align}
(u)^\kappa (\nabla u)^\nu = u^{\kappa_1} \bar{u}^{\kappa_2}
u^{\nu_1}_{x_1} \bar{u}^{\nu_2}_{x_1}...u^{\nu_{2n-1}}_{x_n}
\bar{u}^{\nu_{2n}}_{x_n}, \label{pol-not}
\end{align}
$|\kappa|=\kappa_1+\kappa_2$, $|\nu|=\nu_1+...+\nu_{2n}$. By Lemma
\ref{NLE-Besov}, we have
\begin{align}
& \|(u)^\kappa (\nabla u)^\nu\| _{\ell^{1,s-1/2}_\triangle
\ell^1_\alpha (L^2_{t,x}(\mathbb{R} \times Q_\alpha))}
\nonumber\\
& \lesssim \sum_{i=0,1} \|\nabla^i u\|_{\ell^{1,s-1/2}_\triangle
\ell^\infty_\alpha (L^2_{t,x}(\mathbb{R} \times Q_\alpha))}
\sum_{k=0,1} \|\nabla^k u\|^{|\kappa|+|\nu|-1}_{\ell^{1}_\triangle
\ell^{|\kappa|+|\nu|-1}_\alpha (L^\infty_{t,x}(\mathbb{R} \times
Q_\alpha))} . \label{proof-thm-2}
\end{align}
Hence, if $u\in \mathscr{D}$, in view of \eqref{proof-thm-1} and
\eqref{proof-thm-2}, we have
\begin{align}
& \lambda_1 ( \mathscr{T} u)  \lesssim  \|u_0\|_{B^{s}_{2,1}} +
\sum_{m+1\le |\kappa|+|\nu| \le M+1} \varrho ^{|\kappa|+|\nu|}.
\label{proof-thm-3}
\end{align}
In view of the estimate for the maximal function as in Proposition
\ref{prop-max-2}, one has that
\begin{align}
& \lambda_2 ( S(t) u_0)  \lesssim  \|u_0\|_{B^{s}_{2,1}}.
\label{proof-thm-4}
\end{align}
and for $i=0,1$,
\begin{align}
\|\mathscr{A} \nabla^i F \|_{\ell^{1}_\triangle \ell^{2^*}_\alpha
(L^\infty_{t,x}(\mathbb{R} \times Q_\alpha))}  & \le
\sum^\infty_{j=0} \int_{\mathbb{R}} \|S(t-\tau)(\triangle_j \nabla^i
F)(\tau)\|_{ \ell^{2^*}_\alpha (L^\infty_{t,x}(\mathbb{R}
\times Q_\alpha))} \nonumber\\
& \lesssim \sum^\infty_{j=0} \int_{\mathbb{R}} \|(\triangle_j
\nabla^i
F)(\tau)\|_{H^{s-3/2}(\mathbb{R}^n) } d\tau \nonumber\\
& \lesssim \sum^\infty_{j=0} 2^{(s-1/2)j} \int_{\mathbb{R}}
\|(\triangle_j F)(\tau)\|_{L^2(\mathbb{R}^n)} d\tau.
\label{proof-thm-5}
\end{align}
Hence, by \eqref{proof-thm-4} and \eqref{proof-thm-5},
\begin{align}
& \lambda_2 ( \mathscr{T} u)  \lesssim  \|u_0\|_{B^{s}_{2,1}} +\|F\|
_{\ell^{1,s-1/2}_\triangle \ell^1_\alpha (L^1_t L^2_{x}(\mathbb{R}
\times Q_\alpha))}.  \label{proof-thm-6}
\end{align}
Similar to \eqref{proof-thm-6}, in view of Proposition
\ref{prop-max-4}, we have
\begin{align}
& \lambda_3 ( \mathscr{T} u)  \lesssim  \|u_0\|_{B^{s}_{2,1}} +\|F\|
_{\ell^{1,s-1/2}_\triangle \ell^1_\alpha (L^1_t L^2_{x}(\mathbb{R}
\times Q_\alpha))}.  \label{proof-thm-7}
\end{align}
In view of Lemma \ref{NLE-Besov}, we have
\begin{align}
 \|(u)^\kappa (\nabla u)^\nu\| _{\ell^{1,s-1/2}_\triangle
\ell^1_\alpha (L^1_t L^2_{x}(\mathbb{R} \times Q_\alpha))}  &
\lesssim \sum_{k=0,1} \|\nabla^k
u\|^{|\kappa|+|\nu|-1}_{\ell^{1}_\triangle
\ell^{(|\kappa|+|\nu|-1)}_\alpha (L^{2(|\kappa|+|\nu|-1)}_t
L^\infty_{x}(\mathbb{R} \times Q_\alpha))}
\nonumber\\
& \ \ \ \ \times \sum_{i=0,1} \|\nabla^i
u\|_{\ell^{1,s-1/2}_\triangle \ell^\infty_\alpha
(L^2_{t,x}(\mathbb{R} \times Q_\alpha))}. \label{proof-thm-8}
\end{align}
Hence, if $u \in \mathscr{D}$, we have
\begin{align}
& \lambda_2 ( \mathscr{T} u) + \lambda_3 ( \mathscr{T} u) \lesssim
\|u_0\|_{B^{s}_{2,1}} + \sum_{m+1\le |\kappa|+|\nu| \le M+1} \varrho
^{|\kappa|+|\nu|}. \label{proof-thm-9}
\end{align}
Repeating the procedures as in the above, we obtain that there
exists $u\in \mathscr{D}$ satisfying the integral equation
$\mathscr{T} u=u$, which finishes the proof of Theorem \ref{GWP-nD}.
$\hfill \Box$

\section{Proof of Theorem \ref{LWP-Mod}} \label{proof-lwp-1d}
We begin with the following
\begin{lem} \label{lem-mod-1}
Let $\mathscr{A}$ be as in \eqref{S(t)}. There exists a constant
$C(T)>1$ which depends only on $T$ and $n$ such that
\begin{align}
\sum_{i=0,1}\|\mathscr{A} \nabla^i F\|_{\ell^{1}_\Box \ell^2_\alpha
(L^\infty_{t,x}([0,T]\times Q_\alpha))} \le C(T)   \|
F\|_{\ell^{1,3/2}_\Box \ell^1_\alpha (L^1_t L^2_{x}([0,T]\times
Q_\alpha))}. \label{modmax1}
\end{align}
\end{lem}
{\bf Proof.} Using Minkowski's inequality and Proposition
\ref{prop-max-1},
\begin{align}
& \|\mathscr{A} \nabla^i F\|_{\ell^1_\Box \ell^2_\alpha
(L^\infty_{t,x}([0,T]\times Q_\alpha))} \nonumber\\
& \le \sum_{k\in \mathbb{Z}^n} \left( \sum_{\alpha\in \mathbb{Z}^n}
\left( \int^T_0 \|S(t-\tau)\Box_k \nabla^i F
(\tau)\|_{L^\infty_{t,x}([0,T]\times Q_\alpha)} d\tau \right)^2
\right)^{1/2} \nonumber\\
& \le \sum_{k\in \mathbb{Z}^n} \int^T_0 \left( \sum_{\alpha\in
\mathbb{Z}^n}   \|S(t-\tau)\Box_k \nabla^i F
(\tau)\|^2_{L^\infty_{t,x}([0,T]\times Q_\alpha)}  \right)^{1/2}
d\tau \nonumber\\
& \le \sum_{k\in \mathbb{Z}^n} \int^T_0    \|\Box_k \nabla^i F
(\tau)\|_{M^{1/2}_{2,1}}  d\tau. \label{modmax2}
\end{align}
It is easy to see that for $i=0,1$,
\begin{align}
\|\Box_k \nabla^i F \|_{M^{1/2}_{2,1}}  & \lesssim \langle
k\rangle^{3/2} \|\Box_k  F \|_{L^2(\mathbb{R}^n)}  \le \langle
k\rangle^{3/2} \sum_{\alpha\in \mathbb{Z}^n}\|\Box_k  F
\|_{L^2(Q_\alpha)} . \label{modmax3}
\end{align}
By \eqref{modmax2} and \eqref{modmax3}, we immediately have
\eqref{modmax1}.  $\hfill \Box$

\begin{lem} \label{lem-mod-2}
Let $\mathscr{A}$ be as in \eqref{S(t)}. Let $n\ge 2$, $s>0$. Then
we have
\begin{align}
\sum_{i=0,1}\|\nabla^i \mathscr{A} F\|_{\ell^{1,s}_\Box
\ell^\infty_\alpha (L^2_{t,x}([0,T]\times Q_\alpha))}  \le \langle
T\rangle^{1/2} \| F\|_{\ell^{1,s}_\Box \ell^1_\alpha
(L^2_{t,x}([0,T]\times Q_\alpha))}. \label{modmax-1}
\end{align}
\end{lem}
{\bf Proof.} In view of Proposition \ref{local-sm-3},  we have
\begin{align}
\|\nabla \mathscr{A} F\|_{\ell^{1,s}_\Box \ell^\infty_\alpha
(L^2_{t,x}([0,T]\times Q_\alpha))} \lesssim  \| F\|_{\ell^{1,s}_\Box
\ell^1_\alpha (L^2_{t,x}([0,T]\times Q_\alpha))}. \label{modmax-2}
\end{align}
By Propositions \ref{local-Str} and \ref{local-Str-2D},
\begin{align}
\|\mathscr{A} F\|_{\ell^1_\Box \ell^\infty_\alpha
(L^2_{t,x}([0,T]\times Q_\alpha))} & \lesssim \|
F\|_{\ell^{1,s}_\Box \ell^1_\alpha (L^1_tL^2_{x}([0,T]\times
Q_\alpha))}
\nonumber\\
& \le  T^{1/2} \| F\|_{\ell^{1,s}_\Box \ell^1_\alpha
(L^2_{t,x}([0,T]\times Q_\alpha))}.
 \label{modmax-3}
\end{align}
By \eqref{modmax-2} and \eqref{modmax-3} we immediately have
\eqref{modmax-1}. $\hfill \Box$

\begin{lem} \label{lem-mod-3}
Let $n\ge 2$, $S(t)$ be as in \eqref{S(t)}.  Then we have for $i=0,
1$,
\begin{align}
\|\nabla^i S(t) u_0\|_{\ell^{1,s}_\Box \ell^\infty_\alpha
(L^2_{t,x}([0,T]\times Q_\alpha))} & \lesssim
\|u_0\|_{M^{s+1/2}_{2,1}}, \quad n\ge 3, \label{modmax-S-1}\\
\|\nabla^i S(t) u_0\|_{\ell^{1,s}_\Box \ell^\infty_\alpha
(L^2_{t,x}([0,T]\times Q_\alpha))} & \lesssim
\|u_0\|_{M^{s+1/2}_{2,1} \cap \dot H^{-1/2}}, \quad n=2.
\label{modmax-S-2}
\end{align}
\end{lem}
{\bf Proof.} \eqref{modmax-S-1} follows from Corollary
\ref{local-sm-2}. For $n=2$,  by Proposition \ref{local-Str-2DH}, we
have the result, as desired. $\hfill \Box$

\begin{lem} \label{NLE-modulation}
Let $n\ge 2$, $s>0$, $L\in \mathbb{N}, L\ge 3$. Let $1\le p, p_i, q,
q_i \le \infty$ satisfy $1/p = 1/p_1+(L-1)/p_2$ and $1/q=1/q_1+
(L-1)/q_2$. We have
\begin{align}
 \|v_1...v_L\|_{\ell^{1,s}_\Box \ell^1_\alpha (L^q_t L^p_x
(I\times Q_\alpha))}  \lesssim &
 \sum^L_{l=1}  \|v_l\|_{\ell^{1,s}_\Box \ell^\infty_\alpha (L^{q_1}_t L^{p_1}_x
(I \times Q_\alpha))} \nonumber\\
&  \  \times \prod_{i\not= l, \, i=1,...,L} \|v_i\|_{\ell^{1}_\Box
\ell^{L-1}_\alpha (L^{q_2}_t L^{p_2}_x (I\times Q_\alpha))}.
\label{NLE-mod1}
\end{align}
\end{lem}
{\bf Proof.} Using the identity
\begin{align}
 v_1...v_L = \sum_{k_1,..., k_L \in \mathbb{Z}^n} \Box_{k_1} v_1...  \Box_{k_L} v_L
\label{NLE-mod2}
\end{align}
and noticing the fact that
\begin{align}
\Box_k( \Box_{k_1} v_1...  \Box_{k_L} v_L) =0, \ \  |k-k_1-...-k_L|
\ge C(L,n), \label{NLE-mod3}
\end{align}
we have
\begin{align}
& \|v_1...v_L\|_{\ell^{1,s}_\Box \ell^1_\alpha (L^q_t L^p_x (I\times
Q_\alpha))}   \nonumber\\
& = \sum_{k\in \mathbb{Z}^n} \langle k\rangle^{s} \sum_{\alpha \in
\mathbb{Z}^n} \|\Box_k( v_1... v_L)\|_{L^q_t
L^p_x (I\times Q_\alpha)}  \nonumber\\
& \le \sum_{k_1,...,k_n \in \mathbb{Z}^n} \sum_{|k-k_1-...-k_L|\le
C} \langle k\rangle^{s} \sum_{\alpha \in \mathbb{Z}^n} \|\Box_k(
\Box_{k_1} v_1... \Box_{k_L} v_L)\|_{L^q_t L^p_x (I\times
Q_\alpha)}. \label{NLE-mod4}
\end{align}
Similar to \eqref{NLE-7} and noticing the fact that
$\|\mathscr{F}^{-1}\sigma_k\|_{L^1{(\mathbb{R}^n})} \lesssim 1$, we
have
\begin{align}
 \sum_{\alpha\in \mathbb{Z}^n} \|\,\Box_k
f\|_{L^q_tL^p_{x}(I\times Q_\alpha)} & =\sum_{\alpha\in
\mathbb{Z}^n} \|(\mathscr{F}^{-1}\sigma_k)*
f\|_{L^q_tL^p_{x}(I\times Q_\alpha)}
 \nonumber\\
& \le \int_{\mathbb{R}^n}
|(\mathscr{F}^{-1}\sigma_k)(y)|\left(\sum_{\alpha\in \mathbb{Z}^n}
\| f(t, x-y)\|_{L^q_tL^p_{x}(I\times Q_\alpha)}\right)dy
 \nonumber\\
& \le \sup_{y\in \mathbb{R}^n} \sum_{\alpha\in \mathbb{Z}^n} \| f(t,
x-y)\|_{L^q_tL^p_{x}(I\times Q_\alpha)} \|\mathscr{F}^{-1}\sigma_k
\|_{L^1(\mathbb{R}^n)}
 \nonumber\\
& \lesssim \sum_{\alpha\in \mathbb{Z}^n} \|
f\|_{L^q_tL^p_{x}(I\times Q_\alpha)}. \label{sum-invariance}
\end{align}
By \eqref{NLE-mod4} and  \eqref{sum-invariance}, we have
\begin{align}
& \|v_1...v_L\|_{\ell^{1,s}_\Box \ell^1_\alpha (L^q_t L^p_x (I\times
Q_\alpha))} \nonumber\\
& \le \sum_{k_1,...,k_n \in \mathbb{Z}^n} \sum_{|k-k_1-...-k_L| \le
C} \langle k \rangle^{s} \sum_{\alpha \in \mathbb{Z}^n} \|\Box_{k_1}
v_1...
\Box_{k_L} v_L\|_{L^q_t L^p_x (I\times Q_\alpha)} \nonumber\\
& \lesssim \sum_{k_1,...,k_n \in \mathbb{Z}^n}  (\langle k_1
\rangle^s +...+ \langle k_L \rangle^{s}) \sum_{\alpha \in
\mathbb{Z}^n} \|\Box_{k_1} v_1... \Box_{k_L} v_L\|_{L^q_t L^p_x
(I\times Q_\alpha)}. \label{NLE-mod5}
\end{align}
By H\"older's inequality,
\begin{align}
& \sum_{k_1,...,k_n \in \mathbb{Z}^n}  \langle k_1 \rangle^s
\sum_{\alpha \in \mathbb{Z}^n} \|\Box_{k_1} v_1... \Box_{k_L}
v_L\|_{L^q_t L^p_x (I\times Q_\alpha)} \nonumber\\
& \le \sum_{k_1,...,k_n \in \mathbb{Z}^n}  \langle k_1 \rangle^s
\sum_{\alpha \in \mathbb{Z}^n} \|\Box_{k_1} v_1\|_{L^{q_1}_t
L^{p_2}_x (I\times Q_\alpha)} \prod^L_{i=2} \|\Box_{k_i}
v_i\|_{L^{q_2}_t L^{p_2}_x (I\times Q_\alpha)} \nonumber\\
& \le \|v_1\|_{\ell^{1,s}_\Box \ell^\infty_\alpha (L^{q_1}_t
L^{p_2}_x (I\times Q_\alpha))} \sum_{k_2,...,k_n \in \mathbb{Z}^n}
\sum_{\alpha \in \mathbb{Z}^n}  \prod^L_{i=2} \|\Box_{k_i}
v_i\|_{L^{q_2}_t L^{p_2}_x (I\times Q_\alpha)} \nonumber\\
&
\le \|v_1\|_{\ell^{1,s}_\Box \ell^\infty_\alpha (L^{q_1}_t
L^{p_2}_x (I\times Q_\alpha))}
\sum_{k_2,...,k_n \in \mathbb{Z}^n}
 \prod^L_{i=2} \|\Box_{k_i}
v_i\|_{\ell^{L-1}_\alpha (L^{q_2}_t L^{p_2}_x (I\times Q_\alpha))}
\nonumber\\
& \le \|v_1\|_{\ell^{1,s}_\Box \ell^\infty_\alpha (L^{q_1}_t
L^{p_2}_x (I\times Q_\alpha))}
 \prod^L_{i=2} \|
v_i\|_{\ell^1_\Box \ell^{L-1}_\alpha (L^{q_2}_t L^{p_2}_x (I\times
Q_\alpha))}. \label{NLE-mod6}
\end{align}
The result follows. $\hfill \Box$

\noindent {\bf Proof of Theorem \ref{LWP-Mod}.} Denote
\begin{align}
\lambda_1(v)& = \sum_{i=0,1} \|\nabla^i v\|_{\ell^{1,3/2}_\Box
\ell^\infty_\alpha (L^2_{t,x}([0,T]\times Q_\alpha))}, \nonumber\\
\lambda_2(v) & = \sum_{i=0,1} \|\nabla^i v\|_{\ell^1_\Box
\ell^2_\alpha (L^\infty_{t,x}([0,T]\times Q_\alpha))}. \nonumber
\end{align}
 Put
\begin{align}
& \mathscr{D} = \left\{ u: \;   \lambda_1 ( u) +  \lambda_2 ( u) \le
\varrho \right\}. \label{metricsp-mod}
\end{align}
Let $\mathscr{T}$ be as in \eqref{map}. We will show that
$\mathscr{T}: \mathscr{D} \to \mathscr{D}$ is a contraction mapping.
First, we consider the case $n\ge 3$. Let $u\in \mathscr{D}$. By
Lemmas \ref{lem-mod-2} and \ref{lem-mod-3}, we have
\begin{align}
\lambda_1( \mathscr{T} u)  \lesssim \|u_0\|_{M^{2}_{2,1}} + \langle
T\rangle^{1/2} \| F\|_{\ell^{1,3/2}_\Box \ell^1_\alpha
(L^2_{t,x}([0,T]\times Q_\alpha))}. \label{LWP-1}
\end{align}
We use the same notation as in \eqref{pol-not}. We have from Lemma
\ref{NLE-modulation} that
\begin{align}
 \|(u)^\kappa (\nabla u)^\nu \|_{\ell^{1,3/2}_\Box \ell^1_\alpha (L^{2}_{t,x}
([0,T]\times Q_\alpha))}  \lesssim &
  \sum_{i=0,1} \|\nabla^i u\|_{\ell^{1,3/2}_\Box \ell^\infty_\alpha (L^{2}_{t,x}
([0,T] \times Q_\alpha))} \nonumber\\
&  \  \times \sum_{k=0,1} \| \nabla^k
u\|^{|\kappa|+|\nu|-1}_{\ell^{1}_\Box \ell^{|\kappa|+|\nu|-1}_\alpha
(L^{\infty}_{t,x} ([0,T]\times Q_\alpha))}\nonumber\\
&  \lesssim \lambda_1(u) \lambda_2(u)^{|\kappa|+|\nu|-1} \le
\varrho^{|\kappa|+|\nu|}. \label{LWP-2}
\end{align}
Hence, for $n\ge 3$,
\begin{align}
 \lambda_1( \mathscr{T} u)  \lesssim
\|u_0\|_{M^{2}_{2,1}} +  \sum^M_{|\kappa|+|\nu|=m+1}
\varrho^{|\kappa|+|\nu|}. \label{3D-estim}
\end{align}
Next, we consider the estimate of $\lambda_2 (\mathscr{T} u)$. By
Lemma \ref{lem-mod-1} and Proposition \ref{prop-max-1},
\begin{align}
\lambda_2( \mathscr{T} u) & \lesssim \|u_0\|_{M^{3/2}_{2,1}} + C(T)
\| F\|_{\ell^{1,3/2}_\Box \ell^1_\alpha (L^1_t L^2_{x}([0,T]\times
Q_\alpha))}, \label{LWP-3}
\end{align}
which reduces to the estimates of $\lambda_1(\cdot)$ as in
\eqref{LWP-1}. Similarly, for $n=2$,
\begin{align}
 \lambda_1( \mathscr{T} u) + \lambda_2( \mathscr{T} u) \lesssim
\|u_0\|_{M^{2}_{2,1} \cap \dot H^{-1/2}} +
\sum^M_{|\kappa|+|\nu|=m+1} \varrho^{|\kappa|+|\nu|}.
\label{2D-estim}
\end{align}
Repeating the procedures as in the proof of Theorem \ref{GWP-nD}, we
can show our results, as desired. $\hfill \Box$

\section{Proof of Theorem \ref{GWP-Mod}} \label{proof-GWP-Mod}
The proof of Theorem \ref{GWP-Mod} follows an analogous way as that
in Theorems \ref{GWP-nD} and \ref{LWP-Mod} and will be sketched. Put
\begin{align}
& \lambda_1(v)= \sum_{i=0,1} \|\nabla^i v\|_{\ell^{1,s-1/2}_\Box
\ell^\infty_\alpha
(L^2_{t,x}(\mathbb{R} \times Q_\alpha))}, \nonumber\\
& \lambda_2(v)= \sum_{i=0,1} \|\nabla^i v\|_{\ell^1_\Box
\ell^m_\alpha
(L^\infty_{t,x}(\mathbb{R} \times Q_\alpha))}, \nonumber\\
& \lambda_3(v)= \sum_{i=0,1} \|\nabla^i v\|_{\ell^1_\Box
\ell^m_\alpha (L^{2m}_t L^\infty_{x}(\mathbb{R} \times Q_\alpha))}.
 \nonumber
\end{align}
Put
\begin{align}
& \mathscr{D}= \left\{ u: \;   \lambda_1 ( u)+\lambda_2 ( u)+
\lambda_3 ( u) \le \varrho \right\}. \label{metricsp-gMod}
\end{align}
Let $\mathscr{T}$ be as in \eqref{map}. We show that $\mathscr{T}:
\mathscr{D} \to \mathscr{D}$. We only consider the case $n\ge 3$.
 It follows from Lemma \ref{lem-mod-3} and \ref{LE-Besov-1}
that
\begin{align}
 \lambda_1(\mathscr{T} u)  \lesssim \|u_0\|_{M^{s}_{2,1}} +
\|F\|_{\ell^{1,s-1/2}_\Box \ell^1_\alpha (L^2_{t,x}(\mathbb{R}
\times Q_\alpha))}. \label{lambdamod-1}
\end{align}
Using Lemma \ref{NLE-modulation} and similar to \eqref{LWP-2}, one
sees that if $u\in \mathscr{D}$, then
\begin{align}
 \lambda_1(\mathscr{T} u)  \lesssim \|u_0\|_{M^{s}_{2,1}} +
 \sum_{m+1\le |\kappa|+|\nu| \le M+1} \varrho^{|\kappa|+|\nu|}.
\label{lambdamod-6}
\end{align}
Using Proposition \ref{prop-max-5} and combining the proof of
\eqref{proof-thm-5}--\eqref{proof-thm-7}, we see that
\begin{align}
 \lambda_2(\mathscr{T} u) +\lambda_3(\mathscr{T} u)  \lesssim \|u_0\|_{M^{s}_{2,1}} +
 \sum_{m+1\le |\kappa|+|\nu| \le M+1} \varrho^{|\kappa|+|\nu|}.
\label{lambdamod-7}
\end{align}
The left part of the proof is analogous to  that of Theorems
\ref{GWP-nD} and \ref{LWP-Mod} and the details are omitted.
$\hfill\Box$

\section{Proof of Theorem \ref{GWP-1D}} \label{proof-gwp-1D}

We prove Theorem \ref{GWP-1D} by following some ideas as in Molinet
and Ribaud \cite{MR} and Wang and Huang \cite{WaHu}. The following
is the estimates for the solutions of the linear Schr\"odinger
equation, see  \cite{KePoVe,MR,WaHu}. Recall that $\triangle_{j}:=
\mathscr{F}^{-1} \delta(2^{-j}\,\cdot) \mathscr{F}$, $j\in
\mathbb{Z}$ and $\delta (\cdot )$ is as in Section
\ref{functionspace}.

\begin{lem} \label{cor4}
 Let
$g\in \mathscr{S}(\mathbb{R})$,$ f\in\mathscr{S}(\mathbb{R}^{2})$,
$4\le p<\infty$. Then we have
\begin{align}
\|\triangle_{j} S(t)g\|_{L_{t}^{\infty}L_{x}^{2} \, \cap \,
L^6_{x,t}} & \lesssim
\|\triangle_{j}g\|_{L^{2}},\label{cor-1}\\
 \|\triangle_{j}S(t)g\|_{L_x^p L_t^\infty} &
\lesssim 2^{j(\frac{1}{2}-\frac{1}{p})}
\|\triangle_{j}g\|_{L^{2}},\label{cor-4}\\
\|\triangle_{j} S(t)g\|_{L_{x}^{\infty}L_{t}^{2}} & \lesssim
2^{-j/2}\|\triangle_{j}g\|_{L^{2}},\label{cor-7}
\end{align}
\begin{align}
\|\triangle_{j} \mathscr{A}f \|_{L_{t}^{\infty}L_{x}^{2} \, \cap \,
L^6_{x,t} } & \lesssim
\|\triangle_{j}f\|_{L^{6/5}_{x,t}}, \label{cor-2}\\
\|\triangle_{j} \mathscr{A}f \|_{L_x^p L_t^\infty } & \lesssim
2^{j(\frac{1}{2}-\frac{1}{p})}
\|\triangle_{j}f\|_{L^{6/5}_{x,t}}, \label{cor-5}\\
\|\triangle_{j} \mathscr{A}f\|_{L_{x}^{\infty}L_{t}^{2}} & \lesssim
2^{-j/2} \|\triangle_{j}f\|_{L^{6/5}_{x,t}}, \label{cor-8}
\end{align}
and
 \begin{align}
\|\triangle_{j}
\mathscr{A}(\partial_{x}f)\|_{L_{t}^{\infty}L_{x}^{2} \, \cap \,
L^6_{x,t} } & \lesssim 2^{j/2}\|\triangle_{j}f\|_{L^1_x
L^2_t}, \label{cor-3}\\
 \|\triangle_{j}
\mathscr{A}(\partial_{x}f)\|_{L_x^p L_t^\infty  } & \lesssim 2^{j/2}
2^{j(\frac{1}{2}-\frac{1}{p})}
\|\triangle_{j}f\|_{L^1_x L^2_t}, \label{cor-6}\\
\|\triangle_{j}
\mathscr{A}(\partial_{x}f)\|_{L_{x}^{\infty}L_{t}^{2}} & \lesssim
\|\triangle_{j}f\|_{L_{x}^{1}L_{t}^{2}}.\label{cor-9}
\end{align}
\end{lem}

For convenience, we write for any Banach function space $X$,
$$
\|f\|_{\ell_\triangle^{1,s}  (X)}= \sum_{j\in \mathbb{Z}}
2^{js}\|\triangle_j f\|_{X}, \quad \|f\|_{\ell_\triangle^1  (X)}:=
\|f\|_{\ell_\triangle^{1,0} (X)}.
$$

\begin{lem}\label{nonlinear-estimate-1}
Let $s> 0$,  $1\le p, p_i, \gamma, \gamma_i \le \infty$ satisfy
\begin{align}
\frac{1}{p}= \frac{1}{p_1}+...+\frac{1}{p_N}, \quad
\frac{1}{\gamma}= \frac{1}{\gamma_1}+...+   \frac{1}{\gamma_N}.
\label{p-gamma}
\end{align}
Then
\begin{align}
\left\|u_1 ... u_N \right\|_{\ell_\triangle^{1,s} (L^p_x
L^{\gamma}_t )} & \lesssim \|u_1 \|_{\ell_\triangle^{1,s} (L^{p_1}_x
L^{\gamma_1}_t )} \prod^N_{i=2} \| u_i\|_{\ell_\triangle^{1}
(L^{p_i}_x L^{\gamma_i}_t )}
\nonumber\\
 & \quad +\|u_2 \|_{\ell_\triangle^{1,s} ( L^{p_2}_x L^{\gamma_2}_t)}
\prod_{i\not= 2,\, i=1,...,N} \| u_i\|_{\ell_\triangle^{1} (L^{p_i}_x L^{\gamma_i}_t )} \nonumber\\
& \quad   +... +  \prod^{N-1}_{i=1} \|u_i\|_{\ell_\triangle^{1} (
L^{p_i}_x L^{\gamma_i}_t)}   \| u_N \|_{\ell_\triangle^{1,s}
(L^{p_N}_x L^{\gamma_N}_t )}, \label{pf-7-1}
\end{align}
and in particular, if $u_1=...=u_N=u$, then
\begin{align}
\left\| u^N \right\|_{\ell_\triangle^{1,s} ( L^p_xL^{\gamma}_t)} &
\lesssim \|u \|_{\ell_\triangle^{1,s} ( L^{p_1}_xL^{\gamma_1}_t)}
\prod^N_{i=2} \| u\|_{\ell_\triangle^{1} (
L^{p_i}_xL^{\gamma_i}_t)}. \label{pf-7-2}
\end{align}
Substituting the spaces $ L^{p}_xL^{\gamma}_t$ and $L^{p_i}_x
L^{\gamma_i}_t $ by $L^{\gamma}_tL^{p}_x  $ and $L^{\gamma_i}_t
L^{p_i}_x $, respectively, \eqref{pf-7-1} and \eqref{pf-7-2} also
holds.
\end{lem}

\noindent {\bf Proof.} We only consider the case $N=2$ and the case
$N>2$ can be handled in a similar way. We have
\begin{align}
u_1u_2 &= \sum^\infty_{r=-\infty} [(S_{r+1}u_1) (S_{r+1}u_2)-
(S_{r}u_1) (S_{r}u_2)] \nonumber\\
&= \sum^\infty_{r=-\infty} [(\triangle_{r+1} u_1) (S_{r+1}u_2)+
(S_{r}u_1) (\triangle_{r+1}u_2)], \label{pf-7-3}
\end{align}
and
\begin{align}
\triangle_j(u_1u_2) = \triangle_j \Big(\sum_{r\ge j-10}
[(\triangle_{r+1} u_1) (S_{r+1}u_2)+ (S_{r}u_1)
(\triangle_{r+1}u_2)]\Big). \label{pf-7-4}
\end{align}
We may assume, without loss of generality that there is only the
first term in the right hand side of \eqref{pf-7-4} and the second
term can be handled in the same way. It follows from Bernstein's
estimate, H\"older's  and Young's inequalities that
\begin{align}
\sum_{j\in \mathbb{Z}} 2^{sj}
\|\triangle_j(u_1u_2)\|_{L^p_xL^\gamma_t} & \lesssim \sum_{j\in
\mathbb{Z}} 2^{sj} \sum_{r\ge j-10} \| (\triangle_{r+1} u_1)
(S_{r+1}u_2) \|_{L^p_xL^\gamma_t} \nonumber\\
& \lesssim \sum_{j\in \mathbb{Z}} 2^{sj} \sum_{r\ge j-10} \|
\triangle_{r+1} u_1\|_{L^{p_1}_xL^{\gamma_1}_t}
\|S_{r+1}u_2 \|_{L^{p_2}_xL^{\gamma_2}_t} \nonumber\\
& \lesssim \sum_{j\in \mathbb{Z}} 2^{s(j-r)} \sum_{r\ge j-10} 2^{rs}
\| \triangle_{r+1} u_1\|_{L^{p_1}_xL^{\gamma_1}_t}
\|S_{r+1}u_2 \|_{L^{p_2}_xL^{\gamma_2}_t} \nonumber\\
& \lesssim \| u_1\|_{\ell_\triangle^{1,s}(L^{p_1}_xL^{\gamma_1}_t)}
\|u_2 \|_{\ell_\triangle^1(L^{p_2}_xL^{\gamma_2}_t)}, \label{pf-7-5}
\end{align}
which implies the result, as desired. $\quad\quad\Box$
\begin{rem} \label{rem-on-lemma-3.1}
\rm One easily sees that \eqref{pf-7-2} can be slightly improved by
\begin{align}
\left\| u^N \right\|_{\ell^{1,s} (L^p_xL^{\gamma}_t )} & \lesssim
\|u \|_{\ell^{1,s} (L^{p_1}_xL^{\gamma_1}_t )} \prod^N_{i=2} \|
u\|_{ L^{p_i}_xL^{\gamma_i}_t }. \label{pf-7-2-1}
\end{align}
In fact, from Minkowski's inequality it follows that
\begin{align}
\|S_{r}u \|_{L^{p}_xL^{\gamma}_t} \lesssim  \|u
\|_{L^{p}_xL^{\gamma}_t}. \label{pf-7-2-2}
\end{align}
From \eqref{pf-7-5} and \eqref{pf-7-2-2} we get \eqref{pf-7-2-1}.
\end{rem}

\noindent {\bf Proof of Theorem \ref{GWP-1D}.} We can assume,
without loss of generality that
\begin{align}
 F(u, \bar{u}, u_x, \bar{u}_x)=
\sum_{m+1\le \kappa+\nu \le M+1} \lambda_{\kappa\nu} u^\kappa
u^\nu_x \label{simple-poly}
\end{align}
and the general case can be handled in the same way.

{\it Step} 1. We consider the case $m>4.$ Recall that
\begin{align}
 \|u\|_{X} & =\sup_{s_m\le s \le \tilde{s}_{M}} \sum_{i=0,1}
\sum_{j\in \mathbb{Z}}
|\!|\!|\partial^i_x \triangle_j u |\!|\!|_s,   \label{1d.2} \\
|\!|\!|\triangle_j v |\!|\!|_s   :=&   2^{s j} (\|\triangle_j v
\|_{L^\infty_tL^2_x \,\cap\,  L^{6}_{x,t}}   + 2^{j/2}
\|\triangle_j v\|_{L^\infty_x L^2_t})\nonumber\\
& + 2^{(s-\tilde{s}_m) j}\|\triangle_j v\|_{L_{x}^{m}L_{t}^{\infty}}
+2^{(s-\tilde{s}_{M})j} \|\triangle_j v\|_{L_x^{M}L_t^\infty}.
\end{align}
Considering the mapping
\begin{align}
& \mathscr{T}:  u(t) \to  S(t)u_0 -{\rm
i}\mathscr{A}F(u,\bar{u},u_x,\bar{u}_x),
\end{align}
we will show that $\mathscr{T}:X\to X$ is a contraction mapping. We
have
\begin{align}
\| \mathscr{T}  u(t)\|_X \lesssim \| S(t)u_0\|_X +\|\mathscr{A}
F(u,\bar{u},u_x,\bar{u}_x)\|_X. \label{1d.5}
\end{align}
In view of \eqref{cor-1}, \eqref{cor-4} and \eqref{cor-7} we have,
\begin{align}
|\!|\!|\partial^i_x \triangle_j S(t)u_0 |\!|\!|_s  \lesssim 2^{sj}
\|\partial^i_x \triangle_j u_0\|_2.
\end{align}
It follows that
\begin{align}
\| S(t)u_0 \|_X  \lesssim \sup_{s_m\le s\le
\tilde{s}_{M}}\sum_{i=0,1} \sum_{j\in \mathbb{Z}} 2^{sj}
\|\partial^i_x \triangle_j u_0\|_2 \lesssim \|u_0\|_{\dot
B^{s_m}_{2,1} \cap \dot B^{1+\tilde{s}_{M}}_{2,1}}. \label{1d.7}
\end{align}
We now estimate $\|\mathscr{A} F(u,\bar{u},u_x,\bar{u}_x)\|_X.$ We
have from \eqref{cor-2}, \eqref{cor-5} and \eqref{cor-8} that
\begin{align}
|\!|\!| \triangle_j (\mathscr{A} F(u,\bar{u},u_x,\bar{u}_x))
|\!|\!|_s  \lesssim 2^{sj} \|\triangle_j
F(u,\bar{u},u_x,\bar{u}_x)\|_{L^{6/5}_{x,t}}. \label{1d.8}
\end{align}
From \eqref{cor-3}, \eqref{cor-6} and \eqref{cor-9} it follows that
\begin{align}
|\!|\!| \triangle_j (\mathscr{A}\partial_x
F(u,\bar{u},u_x,\bar{u}_x)) |\!|\!|_s  \lesssim 2^{sj}2^{j/2}
\|\triangle_j F(u,\bar{u},u_x,\bar{u}_x)\|_{L^1_x L^2_t}.
\label{1d.9}
\end{align}
Hence, from \eqref{1d.2},  \eqref{1d.8} and \eqref{1d.9} we have
\begin{align}
\|\mathscr{A} F(u,\bar{u},u_x,\bar{u}_x)\|_X \lesssim & \sum_{j\in
\mathbb{Z}} 2^{sj} \|\triangle_j
F(u,\bar{u},u_x,\bar{u}_x)\|_{L^{6/5}_{x,t}} \nonumber\\
& + \sum_{j\in \mathbb{Z}} 2^{sj}2^{j/2} \|\triangle_j
F(u,\bar{u},u_x,\bar{u}_x)\|_{L^1_x L^2_t}=I+II. \label{1d.10}
\end{align}
Now we perform the nonlinear estimates.  By Lemma
\ref{nonlinear-estimate-1},
\begin{align}
I& \lesssim \sum_{m+1\le \kappa+\nu\le M+1} \Big(
\|u\|_{\ell_\triangle^{1,s}(L^6_{x,t})}
\|u\|^{\kappa-1}_{\ell_\triangle^{\,1}(L^{3(\kappa+\nu-1)/2}_{x,t})
}
\|u_x\|^\nu_{\ell_\triangle^{\,1}(L^{3(\kappa+\nu-1)/2}_{x,t}) } \nonumber\\
& \quad \quad\quad \quad \quad\quad
\quad+\|u_x\|_{\ell_\triangle^{1,s}(L^6_{x,t})} \|u_x
\|^{\nu-1}_{\ell_\triangle^{\,1}(L^{3(\kappa+\nu-1)/2}_{x,t}) }
\|u\|^\kappa_{\ell_\triangle^{\,1}(L^{3(\kappa+\nu-1)/2}_{x,t}) }
\Big)\nonumber\\
& \lesssim \sum_{m+1\le \kappa+\nu\le M+1} \big(\sum_{i=0,1}
\|\partial_x^i u\|_{\ell_\triangle^{1,s}(L^6_{x,t})}\big)
\big(\sum_{i=0,1} \|\partial^i_x
u\|^{\kappa+\nu-1}_{\ell_\triangle^{\,1}(L^{3(\kappa+\nu-1)/2}_{x,t})
} \big). \label{1d.11}
\end{align}
For any $m\le \lambda \le M$, we let $\frac{1}{\rho}= \frac{1}{2}-
\frac{4}{3\lambda}$. It is easy to see that the following inclusions
hold:
\begin{align}
L^{\infty}_{t} (\mathbb{R}, \dot H^{s_\lambda}) \cap L^6_t
(\mathbb{R}, \dot H^{s_\lambda}_6) \subset L^{3\lambda/2}_{t}
(\mathbb{R}, \dot H^{s_\lambda}_\rho) \subset L^{3\lambda/2}_{x,t}.
\label{pf-7-7}
\end{align}
More precisely, we have
\begin{align}
\sum_{j\in \mathbb{Z}}\|\triangle_j u\|_{L^{3\lambda/2}_{x,t}} &
\lesssim \sum_{j\in \mathbb{Z}}\|\triangle_j u\|_{
L^{3\lambda/2}_{t}
(\mathbb{R},\, \dot H^{s_\lambda}_\rho)} \nonumber\\
& \lesssim \sum_{j\in \mathbb{Z}}\|\triangle_j u\|^{4/\lambda}_{
L^{6}_{t} (\mathbb{R}, \, \dot H^{s_\lambda}_6)} \|\triangle_j
u\|^{1-4/\lambda}_{
L^{\infty}_{t} (\mathbb{R}, \, \dot H^{s_\lambda})} \nonumber\\
& \lesssim \|u\|^{4/\lambda}_{ \ell^{1,s_\lambda}( L^{6}_{x,t}) } \|
u\|^{1-4/\lambda}_{\ell^{1,s_\lambda}( L^{\infty}_{t} L^2_x )}.
\label{pf-7-8}
\end{align}
Using \eqref{pf-7-8} and noticing that $s_m \le s_{\kappa+\nu-1}\le
s_{M} < \tilde{s}_{M}$, we have
\begin{align}
\|\partial^i_x
u\|^{\kappa+\nu-1}_{\ell_\triangle^{\,1}(L^{3(\kappa+\nu-1)/2}_{x,t})
} & \lesssim  \|\partial^i_x u\|^{4}_{\ell_\triangle^{ 1,
\,s_{\kappa+\nu-1}}(L^{6}_{x,t}) } \|\partial^i_x
u\|^{\kappa+\nu-5}_{\ell_\triangle^{1,\,s_{\kappa+\nu-1}}(L^\infty_t
L^2_x) }
\nonumber\\
& \lesssim \|u\|^{\kappa+\nu-1}_X . \label{1d.12}
\end{align}
Combining \eqref{1d.11} with  \eqref{1d.12}, we have
\begin{align}
I  \lesssim \sum_{m+1\le \kappa+\nu\le M+1} \|u\|^{\kappa+\nu}_X.
\label{1d.13}
\end{align}
Now we estimate $II$. By Lemma \ref{nonlinear-estimate-1},
\begin{align}
II & \lesssim \sum_{m+1\le \kappa+\nu\le M+1} \Big(
\|u\|_{\ell_\triangle^{1,\,s+1/2}(L^\infty_x L^2_t)}
\|u\|^{\kappa-1}_{\ell_\triangle^{\,1}(L^{\kappa+\nu-1}_x
L^\infty_t) }
\|u_x\|^\nu_{\ell_\triangle^{\,1}(L^{\kappa+\nu-1}_x L^\infty_t )} \nonumber\\
& \quad \quad\quad \quad \quad\quad
\quad+\|u_x\|_{\ell_\triangle^{1,\,s+1/2}(L^\infty_x L^2_t)}  \|u_x
\|^{\nu-1}_{L^{\kappa+\nu-1}_x L^\infty_t }
\|u\|^\kappa_{L^{\kappa+\nu-1}_x L^\infty_t }
\Big)\nonumber\\
& \lesssim \sum_{m+1\le \kappa+\nu\le M+1} \big(\sum_{i=0,1}
\|\partial_x^i u\|_{ \ell_\triangle^{1,\,s+1/2}(L^\infty_x L^2_t)
}\big) \big(\sum_{i=0,1} \|\partial^i_x u\|^{\kappa+\nu-1}_{
L^{\kappa+\nu-1}_x L^\infty_t  } \big) \nonumber\\
& \lesssim \sum_{m+1\le \kappa+\nu\le M+1} \|u\|_X \big(\sum_{i=0,1}
\|\partial^i_x u\|^{\kappa+\nu-1}_{ L^{m}_x L^\infty_t \, \cap\,
L^{M}_x L^\infty_t  } \big)
\nonumber\\
& \lesssim \sum_{m+1\le \kappa+\nu\le M+1} \|u\|^{\kappa+\nu}_X.
\label{1d.14}
\end{align}
Collecting \eqref{1d.10},  \eqref{1d.11}, \eqref{1d.13} and
\eqref{1d.14}, we have
\begin{align}
\|\mathscr{A} F(u,\bar{u},u_x,\bar{u}_x)\|_X  & \lesssim
\sum_{m+1\le \kappa+\nu\le M+1} \|u\|^{\kappa+\nu}_X. \label{1d.15}
\end{align}
By \eqref{1d.5}, \eqref{1d.7} and \eqref{1d.15}
\begin{align}
\| \mathscr{T}  u(t)\|_X    \lesssim \|u_0\|_{\dot B^{s_m}_{2,1}
\cap \dot B^{1+\tilde{s}_{M}}_{2,1}}+  \sum_{m+1\le \kappa+\nu\le
M+1} \|u\|^{\kappa+\nu}_X.  \label{1d.16}
\end{align}

{\it Step} 2. We consider the case $m=4$. Recall that
\begin{align}
\|u\|_X & =\sum_{i=0,1} \big(\|\partial^i_x u\|_{L^\infty_tL^2_x
\,\cap\,  L^{6}_{x,t} } + \sup_{s_5\le s \le \tilde{s}_{M}}
\sum_{j\in \mathbb{Z}} |\!|\!|\partial^i_x \triangle_j u |\!|\!|_s
\big). \nonumber
\end{align}
By \eqref{cor-1}, \eqref{cor-4} and \eqref{cor-7},
\begin{align}
\| S(t)u_0 \|_X  \lesssim  \|u_0\|_2+  \sup_{s_5\le s\le
\tilde{s}_{M}}\sum_{i=0,1} \sum_{j\in \mathbb{Z}} 2^{sj}
\|\partial^i_x \triangle_j u_0\|_2 \lesssim \|u_0\|_{
B^{1+\tilde{s}_{M}}_{2,1}}. \label{1d.17}
\end{align}
We now estimate $\|\mathscr{A} F(u,\bar{u},u_x,\bar{u}_x)\|_X.$ By
Strichartz' and H\"older's inequality,  we have
\begin{align}
& \| \mathscr{A} F(u,\bar{u},u_x,\bar{u}_x)\|_{L^\infty_tL^2_x
\,\cap\,  L^{6}_{x,t} }\nonumber\\
 &  \lesssim \sum_{5 \le
\kappa+\nu\le M+1}
\|(|u|+|u_x|)^{\kappa+\nu}\|_{L^{6/5}_{x,t}} \nonumber\\
&  \lesssim \sum_{5 \le \kappa+\nu\le M+1}
\|(|u|+|u_x|)\|_{L^{6}_{x,t}}
\|(|u|+|u_x|)\|^{\kappa+\nu-1}_{L^{3(\kappa+\nu-1)/2}_{x,t}}
\nonumber\\
&  \lesssim \sum_{5 \le \kappa+\nu\le M+1}
\|(|u|+|u_x|)\|_{L^{6}_{x,t}}
\|(|u|+|u_x|)\|^{\kappa+\nu-1}_{L^{6}_{x,t} \cap L^{3M/2}_{x,t} }
\nonumber\\
 & \lesssim \sum_{5 \le \kappa+\nu\le M+1}
\big(\sum_{i=0,1}\|\partial^i_x u \|_{L^{6}_{x,t}}
\big)^{\kappa+\nu} \nonumber\\
& \ \ \ + \sum_{5 \le \kappa+\nu\le M+1}
\big(\sum_{i=0,1}\|\partial^i_x u \|_{L^{6}_{x,t}} \big)
\big(\sum_{i=0,1}\|\partial^i_x u \|_{L^{3M/2}_{x,t}}
\big)^{\kappa+\nu-1}. \label{1d.18}
\end{align}
Applying \eqref{pf-7-8},  we see that \eqref{1d.18} implies that
\begin{align}
& \| \mathscr{A} F(u,\bar{u},u_x,\bar{u}_x)\|_{L^\infty_tL^2_x
\,\cap\,  L^{6}_{x,t} } \lesssim \sum_{5 \le \kappa+\nu\le M+1}
\|u\|_X^{\kappa+\nu}. \label{1d.19}
\end{align}
From Bernstein's estimate and \eqref{cor-3} it follows that
\begin{align}
& \|\partial_x  \mathscr{A} F(u,\bar{u},u_x,\bar{u}_x)
\|_{L^\infty_tL^2_x \,\cap\,  L^{6}_{x,t} } \nonumber\\
&  \le \|P_{\le 1}  (\mathscr{A}\partial_x
F(u,\bar{u},u_x,\bar{u}_x)) \|_{L^\infty_tL^2_x \,\cap\, L^{6}_{x,t}
} \nonumber\\
& \ \ \ + \|P_{>1} (\mathscr{A}\partial_x
F(u,\bar{u},u_x,\bar{u}_x)) \|_{L^\infty_tL^2_x \,\cap\,
L^{6}_{x,t}}
\nonumber\\
& \lesssim \|\mathscr{A} F(u,\bar{u},u_x,\bar{u}_x)
\|_{L^\infty_tL^2_x \,\cap\, L^{6}_{x,t}
} \nonumber\\
& \ \ \ + \sum_{j\gtrsim 1} 2^{j/2} \|\triangle_j
F(u,\bar{u},u_x,\bar{u}_x)\|_{L^1_x L^2_t} \nonumber\\
& \lesssim \|\mathscr{A} F(u,\bar{u},u_x,\bar{u}_x)
\|_{L^\infty_tL^2_x \,\cap\, L^{6}_{x,t}
} \nonumber\\
& \ \ \ + \sum_{j\in \mathbb{Z}} 2^{\tilde{s}_M} 2^{j/2}
\|\triangle_j F(u,\bar{u},u_x,\bar{u}_x)\|_{L^1_x L^2_t} =III+IV.
\label{1d.20}
\end{align}
The estimates of $III$ and $IV$ have been given in \eqref{1d.19} and
\eqref{1d.14}, respectively. We have
\begin{align}
& \|\partial_x  \mathscr{A} F(u,\bar{u},u_x,\bar{u}_x)
\|_{L^\infty_tL^2_x \,\cap\,  L^{6}_{x,t} }   \lesssim \sum_{5 \le
\kappa+\nu\le M+1} \|u\|_X^{\kappa+\nu}.  \label{1d.21}
\end{align}
We have from \eqref{cor-2}--\eqref{cor-8},
\eqref{cor-3}--\eqref{cor-9} that
\begin{align}
\sum_{j\in \mathbb{Z}} |\!|\!| \triangle_j (\mathscr{A}
F(u,\bar{u},u_x,\bar{u}_x)) |\!|\!|_s  &  \lesssim  \sum_{j\in
\mathbb{Z}} 2^{sj} \|\triangle_j
F(u,\bar{u},u_x,\bar{u}_x)\|_{L^{6/5}_{x,t}}, \label{1d.22}\\
\sum_{j\in \mathbb{Z}} |\!|\!| \triangle_j (\mathscr{A}\partial_x
F(u,\bar{u},u_x,\bar{u}_x)) |\!|\!|_s &  \lesssim \sum_{j\in
\mathbb{Z}} 2^{sj}2^{j/2} \|\triangle_j
F(u,\bar{u},u_x,\bar{u}_x)\|_{L^1_x L^2_t} \label{1d.23}
\end{align}
hold for all $s>0$. The right hand side in \eqref{1d.23} has been
estimated by \eqref{1d.14}. So, it suffices to consider the estimate
of the right hand side in \eqref{1d.22}. Let us observe the equality
\begin{align}
F(u,\bar{u},u_x,\bar{u}_x) = \sum_{\kappa+\nu=5} \lambda_{\kappa\nu}
u^\kappa u_x^\nu + \sum_{5<\kappa+\nu \le M+1} \lambda_{\kappa\nu}
u^\kappa u_x^\nu :=V+VI.  \label{1d.24}
\end{align}
For any $s_5\le s\le \tilde{s}_M$,  $VI$ has been handled in
\eqref{1d.11}--\eqref{1d.13}:
\begin{align}
\sum_{5< \kappa+\nu\le M+1}  \sum_{j\in \mathbb{Z}} 2^{sj}
\|\triangle_j (u^\kappa u_x^\nu )\|_{L^{6/5}_{x,t}} & \lesssim
\sum_{5< \kappa+\nu\le M+1} \|u\|^{\kappa+\nu}_X. \label{1d.25}
\end{align}
For the estimate of $V$, we use Remark \ref{rem-on-lemma-3.1}, for
any $s_5\le s\le \tilde{s}_M$,
\begin{align}
\sum_{\kappa+\nu=5}  \sum_{j\in \mathbb{Z}} 2^{sj} \|\triangle_j
(u^\kappa u_x^\nu )\|_{L^{6/5}_{x,t}} & \lesssim
\big(\sum_{i=0,1}\|\partial^i_x u_x\|^4_{L^{6}_{x,t}} \big)
\big(\sum_{i=0,1}\|\partial^i_x u_x\|_{\ell_\triangle^{1,s}
(L^{6}_{x,t})}
\big)  \nonumber\\
& \lesssim \|u\|^5_X.  \label{1d.26}
\end{align}
Summarizing the estimate above,
\begin{align}
\| \mathscr{T}  u(t)\|_X    \lesssim \|u_0\|_{
B^{1+\tilde{s}_{M}}_{2,1}}+  \sum_{5 \le \kappa+\nu\le M+1}
\|u\|^{\kappa+\nu}_X,  \label{1d.27}
\end{align}
whence, we have the results, as desired. $\quad\quad\Box$

\noindent{\bf Acknowledgment.}  This work is supported in part by
the National Science Foundation of China, grants 10571004 and
10621061; and the 973 Project Foundation of China, grant
2006CB805902.


\medskip
\footnotesize
\begin{thebibliography}{31}

\bibitem{BL} J. Bergh and J. L\"{o}fstr\"{o}m,
 Interpolation Spaces,  Springer--Verlag, 1976.

\bibitem{Be-Ta} I. Bejenaru and D. Tataru,  Large data local solutions for the derivative NLS
equation, arXiv:math.AP/0610092 v1.

\bibitem{Bour} J. Bourgain, Fourier transform restriction phenomena for certain
lattice subsets and applications to nonlinear evolution equations,
GAFA, {\bf 3} (1993), 107 - 156 and 209 - 262.

\bibitem{CKSTT} J. Colliander, M. Keel, G. Staffilani, H. Takaoka, and T. Tao,
 A refined global well-posedness result for the
Schr\"odinger equation with derivative, SIAM J. Math. Anal., {\bf
34} (2002), 64--86.

\bibitem{Co-Sa} P. Constantin and J. C. Saut, {\rm Local smoothing properties of
dispersive equations}, J. Amer. Math. Soc., {\bf 1} (1988),
413--446.

\bibitem{Chih1} H. Chihara, Global existence of small solutions to
semilinear Schr\"odinger equations with guage invariance, Publ.
RIMS, {\bf 31} (1995), 731--753.
\bibitem{Chih2} H. Chihara, The initial value problem for cubic
semilinear Schr\"odinger equations with guage invariance, Publ.
RIMS, {\bf 32} (1996), 445--471.
\bibitem{Chih3}
H. Chihara, Gain of regularity for semilinear Schrdinger equations,
Math. Ann. {\bf 315} (1999), 529-567.

\bibitem{Chr} M. Christ, Illposedness of a Schr\"odinger equation
with derivative regularity, Preprint.

\bibitem{Co-Ni1} E. Cordero and F. Nicola. Strichartz
estimates inWiener amalgam spaces for the Schr\"odinger equation.
Math. Nachr., {\bf 281} (2008), 25--41.

\bibitem{Co-Ni2} E Cordero, F Nicola, Metaplectic representation on Wiener amalgam
spaces and applications to the Schr\"odinger equation,  J. Funct.
Anal., {\bf 254} (2008), 506-534.




\bibitem{Fei2} H. G. Feichtinger, Modulation spaces on locally
compact Abelian group, Technical Report, University of Vienna, 1983.
Published in: ``Proc. Internat. Conf. on Wavelet and Applications",
99--140.  New Delhi Allied Publishers, India, 2003.
{http://www.unive.ac.at/nuhag-php/bibtex/
open\_files/fe03-1\_modspa03.pdf}.
\bibitem{Groh} K. Gr\"ochenig,  Foundations of Time-Frequency
Analysis, Birkh\"auser, Boston, MA, 2001.
\bibitem{Ke-Ta} M. Keel and T. Tao, {\rm Endpoint Strichartz estimates,}
Amer. J. Math., {\bf 120} (1998), 955-980.


\bibitem{Grun} A. Gr\"unrock On the Cauchy- and periodic boundary value problem for a
certain class of derivative nonlinear Schr\"odinger equations,
arXiv:math/0006195v1.




\bibitem{KePoVe} C. E.
Kenig, G. Ponce and L. Vega, {\rm Oscillatory integrals and
regularity of dispersive equations,} Indiana Univ. Math. J.,  \bf 40
\rm (1991), 253--288.
\bibitem{KePoVe1} C. E. Kenig, G. Ponce, L. Vega, Small solutions to nonlinear
Schrodinger equation, Ann. Inst. Henri Poincar\'{e}, Sect C, {\bf
10} (1993), 255-288.
\bibitem{KePoVe2} C. E. Kenig, G. Ponce and L. Vega, {\rm Smoothing
effects and local existence theory for the generalized nonlinear
Schr\"odinger equations},  Invent. Math., {\bf 134} (1998),
489--545.

\bibitem{KePoVe3}  C. E. Kenig, G. Ponce, L. Vega,  The Cauchy problem for quasi-linear Schrodinger
equations,  Invent. Math. {\bf 158} (2004), 343--388.

\bibitem{KePoRoVe}  C. E.
Kenig, G. Ponce, C. Rolvent, L. Vega,  The genreal quasilinear
untrahyperbolic Schrodinger equation, Advances in Mathematics {\bf
206} (2006), 402--433.


\bibitem{Klai} S. Klainerman, Long-time behavior of solutions to
nonlinear evolution equations, Arch. Rational Mech. Anal., {\bf 78}
(1982), 73--98.
\bibitem{Kl-Po} S. Klainerman, G. Ponce,  Global small amplitude solutions to
nonlinear evolution equations, Commun. Pure Appl. Math., {\bf 36}
(1983), 133--141.

\bibitem{MR} L. Molinet and F. Ribaud, Well posedness results for
the generalized Benjamin-Ono equation with small initial data, {\rm
J. Math. Pures Appl.,} {\bf 83} (2004),  277-311.

\bibitem{MJ} L. Molinet, J.D.Saut and N.Tzvetkov, Ill-posedness
issues for the Benjamin-Ono equation and related equations, SIAM
J.Math. Anslysis, {\bf 33} (2001), 982-988.

\bibitem{Oz-Zh} T. Ozawa and J. Zhang, Global existence of small classical solutions to nonlinear Schr\"odinger
equations, Ann. I. H. Poincar\'{e}, AN, to appear.
\bibitem{Sjol} P. Sj\"olin, Regularity of solutions to the
Schr\"odinger equations, Duke Math. J., {\bf 55} (1987), 699--715.

\bibitem{SuTo} M. Sugimoto amd N. Tomita, The dilation property of
modulation spaces and their inclusion relation with Besov spaces,
Preprint.

\bibitem{Tao} T. Tao, Spherically averaged endpoint Strichartz estimates for the
two-dimensional Schr\"odinger equation, Commun. PDE, {\bf 25}
(2000), 1471--1485.

\bibitem{Toft} J. Toft, Continuity properties for modulation spaces,
with applications to pseudo-differential calculus, I. J. Funct.
Anal., {\bf 207} (2004), 399--429.

\bibitem{Tr} H. Triebel,  Theory of Function Spaces,
  Birkh\"{a}user--Verlag, 1983.



\bibitem{Wa1} Baoxiang Wang, Lifeng Zhao and Boling Guo, Isometric
decomposition operators, function spaces $E^\lambda_{p,q}$ and
applications to nonlinear evolution equations, J. Funct. Anal., {\bf
233} (2006), 1--39.
\bibitem{WaHe}  Baoxiang Wang and Henryk Hudzik, {\rm The global Cauchy problem for the
NLS and NLKG with small rough data,} J. Differential Equations, {\bf
231} (2007), 36--73.
\bibitem{WaHu} Baoxiang Wang and Chunyan Huang, Frequency-uniform decomposition method for the generalized
BO, KdV and NLS equations,  J. Differential Equations, {\bf 239}
(2007), 213--250.


\bibitem{Vega} L. Vega, The Schr\"odinger equation: pointwise convergence to
the initial data, Proc. Amer. Math. Soc., {\bf 102} (1988),
874--878.



\end{thebibliography}
\end{document}